\documentclass[10pt]{paper}

\usepackage{amssymb,amsmath,enumerate,theorem,epsfig,rotating,graphics,changebar,eepic}%,epsf}%t1enc
\usepackage{amsfonts}
\usepackage{a4,latexsym,parskip}
\usepackage{hyperref}
\hypersetup{pdfauthor=Matthias Sch"utt} \hypersetup{pdftitle=Diss}

\newcommand{\Phoch}[1]{\mathbb{P}^{#1}}
\newcommand{\cohom}[4][]{\mathrm{H}_{#1}^{#2}(#3,#4)}
\newcommand{\C} [1][]{\mathbb{C}^{#1}}
\newcommand{\Q} [1] []{\mathbb{Q}_{#1}}
\newcommand{\N} [1][] {\mathbb{N}_{#1}}

\newcommand{\Z}{\mathbb{Z}}

\newtheorem{Spezial-Theorem}{Theorem}[section]
\newtheorem{Spezial-Proposition}{Proposition}[section]
\theoremstyle{break} 

\ifx\undefined\pdfpageheight
  
\else
  
\fi

\ifx\undefined\pdfpageheight
  
\else
  
\fi

\begin{document}
\setlength{\unitlength}{1cm}

\title{Elliptic fibrations of some extremal K3 surfaces}

\author{Matthias Sch\"utt}

\date{\today}
%\subjclass[2000]{14J27,14J28}
\maketitle

\abstract{This paper is concerned with the construction of
extremal elliptic K3 surfaces. It gives a complete treatment of
those fibrations which can be derived from rational elliptic
surfaces by easy manipulations of their Weierstrass equations. In
particular, this approach enables us to find explicit equations
for 38 semi-stable extremal elliptic K3 fibrations, 32 of which
are indeed defined over $\Q$. They are realized as pull-back of
non-semi-stable extremal rational elliptic surfaces via base
change. This is related to work of J. Top and N. Yui which
exhibited the same procedure for the semi-stable extremal rational
elliptic surfaces. }

\begin{small}

\end{small}

\vspace{0.2cm}

\section{Introduction}
\label{s:intro}

The aim of this paper is to find all extremal elliptic K3
fibrations which can be derived from rational elliptic surfaces by
direct, relatively simple manipulations of the Weierstrass
equations. The main technique will be pull-back by a base change.
We only exclude the general construction involving the induced
$J$-map of the fibration (considered as a base change generally of
degree 24, cf.~\cite[section 2]{MP2}). The base changes we
construct will have degree at most 8. Additionally there is
another effective method if we allow the extremal K3 surface to
have non-reduced fibres. Then we can also manipulate the
Weierstrass equations by adding or transferring common factors,
thus changing the shape of singular fibres rather than introducing
new cusps. In total, this approach will enable us to realize 201
out of the 325 configurations of singular fibres which exist for
extremal elliptic K3 fibrations by the classification of
\cite{SZ}. Note, however, that the configuration does in general
not determine the isomorphism class of the complex surface
(cf.~[1,2,3,10,2*] (No.~148) in Section \ref{s:non-semi-stable}).

For most of this paper, we will concentrate on extremal elliptic
K3 fibrations with only semi-stable fibres. The determination of
the 112 possible configurations of singular fibres originally goes
back to Miranda and Persson \cite{MP2}. For 20 of them,
Weierstrass equations over $\Q$ (or in one case $\Q(\sqrt{5})$)
have been obtained in \cite{I} or \cite{ShCR}, \cite{LY} and
\cite{TY}. These give rise to the elliptic K3 surface with the
maximal singular fibre (the first one from the list in
\cite{MP2}), the Shioda modular ones and those coming from
semi-stable extremal (hence Shioda modular) rational elliptic
surfaces after a quadratic base change. Because of the
construction, it is easy to determine the unique $\C$-isomorphism
classes of these surfaces over $\Q$ using \cite{SZ}.

The main idea of this paper consists in applying a base change
of higher degree to other extremal rational elliptic surfaces
(namely those with three cusps). In the pull-back surface, we
replace the non-reduced singular fibres in such a way by
semi-stable fibres that it is an extremal K3. Indeed, five of the
Shioda modular rational elliptic surfaces can also be obtained in
such a way according to \cite[section 7]{MP1}. Here, we
investigate those base changes which do not factor through the
Shioda modular surfaces. We find Weierstrass equations for 38
further extremal semi-stable elliptic K3 fibrations, only 6 of
which are not defined over $\Q$. Again, the isomorphism classes of
these surfaces can be determined in advance since we know the
Mordell-Weil groups. The surfaces over $\Q$ realize the following
32 configurations of singular fibres in the notation of
\cite{MP2}:

\begin{center}
\begin{tabular}{cccc}
[1,1,1,2,3,16] & [1,1,1,2,5,14] & [1,1,1,3,3,15] & [1,1,1,3,6,12]\cr
[1,1,1,5,6,10] & [1,1,2,2,3,15] & [1,1,2,3,3,14] & [1,1,2,4,4,12]\cr
[1,1,2,4,6,10] & [1,1,3,3,8,8] & [1,1,3,4,6,9] & [1,2,2,2,3,14]\cr
[1,2,2,2,5,12] & [1,2,2,2,7,10] & [1,2,2,3,4,12] & [1,2,2,3,6,10]\cr
[1,2,2,5,6,8] & [1,2,2,6,6,7] & [1,2,3,3,3,12] & [1,2,3,4,4,10]\cr
[1,2,4,4,6,7] & [1,2,4,5,6,6] & [1,3,3,3,5,9] & [1,3,3,5,6,6]\cr
[1,3,4,4,4,8] & [2,2,2,4,6,8] & [2,2,2,3,5,10] & [2,2,3,3,4,10]\cr
[2,2,3,4,5,8] & [2,2,4,4,6,6] & [2,3,3,3,4,9] & [2,3,4,4,5,6]\cr
\end{tabular}
\end{center}

The additional fibrations which result from this approach, but can
only be defined over some quadratic or cubic extension of $\Q$,
are

\begin{center}
\begin{tabular}{ccc}
[1,1,2,2,4,14], & [1,1,2,6,6,8], & [1,2,2,4,5,10]\cr
[1,2,2,4,7,8], & [1,2,3,3,6,9], & [1,2,3,4,6,8].
\end{tabular}
\end{center}

For the non-semi-stable extremal K3 fibrations which can be
derived from rational elliptic surfaces, we produce two long
tables in Section \ref{s:non-semi-stable}.
\vspace{0.1cm}

This paper is organized as follows: After shortly recalling some
basic facts about elliptic surfaces in the next section
(\ref{s:basics}), we will spend the major part of it with
constructing the base changes and giving the resulting equations
for the semi-stable extremal elliptic K3 surfaces (Sections
\ref{s:deflated}, \ref{s:degree 4}, \ref{s:degree 6},
\ref{s:inflating}). Eventually we will also consider the
non-semi-stable fibrations in Section \ref{s:non-semi-stable}
although we will keep their treatment quite concise.

One final remark seems to be in order: There are, of course, many
other ways to produce extremal (or singular) elliptic K3 surfaces.
The perhaps best known is the concept of double sextics as
introduced in \cite{P}. We will not pursue this approach here, so
the interested reader is also referred to \cite{MP3} and
\cite{ATZ} for instructive applications.

\section{Elliptic surfaces over $\mathbf{\Phoch{1}}$ with a section}
\label{s:basics}

An elliptic surface over $\Phoch{1}$, say $Y
\stackrel{r}{\rightarrow} \Phoch{1}$, with a section is given by a
minimal Weierstrass equation
\[
y^2=x^3+Ax+B
\]
where A and B are homogeneous polynomials in the two variables of
$\Phoch{1}$ of degree $4M$ and $6M$, respectively, for some
$M\in\N$. Then the section is the point at $\infty$. The term
minimal refers to the common factors of $A$ and $B$: They are not
allowed to have a common irreducible factor with multiplicity
greater than 3 in $A$ and greater than 5 in $B$. Otherwise we
could cancel these factors by a change of variables. This
convention restricts the singularities of the Weierstrass equation
to rational double points. Then, the surface $Y$ is the minimal
desingularization. In this paper, we are mainly interested in
examples where both $A$ and $B$ have rational coefficients (while
the results of \cite{MP2} will only imply the existence of $A$ and
$B$ over some number field). Of course, we can also assume $A$ and
$B$ to have (minimal) integer coefficients, but we will not go
into detail here.

As announced in the introduction, we are going to pay special
attention to the singular fibres of $Y$. For a general choice of
$A$ and $B$ there will be $12M$ of them (each a rational curve
with a node). The types of the singular fibres, which were first
classified by Kodaira in \cite{K}, can be read off directly from
the $j$-function of $Y$ (cf.~\cite[IV.9, table 4.1]{Si}). Up to a
factor, this is the quotient of $A^3$ by the discriminant $\Delta$
of $Y$ which is defined as $\Delta=-16\,(4A^3+27B^2)$:
\[
j=-\dfrac{1728\,(4A)^3}{\Delta}.
\]
Then $Y$ has singular fibres above the zeroes of $\Delta$ which we
call the \emph{cusps} of $Y$. Let $x_0$ be a cusp and $n$ the
order of vanishing of $\Delta$ at $x_0$. The fibre above $x_0$ is
called \emph{semi-stable} if it is a rational curve with a node
(for $n=1$) or a cycle of $n$ lines, $n>1$. In Kodaira's notation,
this is type $I_n$. The fibre above $x_0$ is semi-stable if and
only if $A$ does not vanish at $x_0$ (i.e.~iff $x_0$ is not a
common zero of $A$ and $B$). On the other hand, we get either a
non-reduced fibre (distinguished by a $^\ast$) over $x_0$ if $A$
and $B$ both vanish at $x_0$ to order at least 2 and 3,
respectively, or an additive fibre of type $II, III$ or $IV$
otherwise.

One common property of the singular fibres is that in every case
the vanishing order of $\Delta$ at the cusp $x_0$ equals the Euler
number of the fibre above $x_0$. Recall that
$\cohom{1}{Y}{\mathcal{O}_Y}=0$ and $p_g=$ dim
$\cohom{2}{Y}{\mathcal{O}_Y}=M-1$, while the canonical divisor
$K_Y=(M-2)~F$ for a general fibre $F$ (cf.~e.g.~\cite[Lecture
III]{M1}). Hence, $Y$ is K3 (resp.~rational) if and only if $M=2$
(resp.~$M=1$). The Euler number of $Y$ equals the sum of the Euler
numbers of its (singular) fibres. By the above considerations,
this coincides with the degree $12M$ of $\Delta$. Hence, we obtain
that $Y$ is a K3 surface if and only if its Euler number equals
24. (On the contrary, $Y$ is rational if and only if $e(Y)=12$.)

Before discussing the effect of a base change on the elliptic
surface $Y$, let us at first introduce the following notation: We
say that a map
\[
\pi: \Phoch{1}\rightarrow\Phoch{1}
\]
has \emph{ramification index} $(n_1,...,n_k)$ at $x_0\in\Phoch{1}$
if $x_0$ has $k$ pre-images under $\pi$ with respective orders
$n_i ~~(i=1,\hdots,k)$. The base change of $Y$ by $\pi$ is simply
defined as the pull-back surface
\[
X \stackrel{\pi\circ
r}{\longrightarrow}\Phoch{1}.
\]
Here, we substitute $\pi$ into the Weierstrass equation and
$j$-function of $Y$ and subsequently normalize to obtain $X$. Let
$\pi$ have ramification index $(n_1,\hdots,n_k)$ at the cusp $x_0$
of $Y$. Then a semi-stable singular fibre of type $I_n$ above
$x_0$ is replaced by $k$ fibres of types
$I_{n_1n},\hdots,I_{n_kn}$ in the pull-back $X$.

For a non-semi-stable singular fibre the substitution %%@
process is non-trivial for two reasons: On the one hand the
Weierstrass equation might simply lose its minimality through the
substitution. On the other hand, the minimalized Weierstrass
equation can still become \emph{inflated}. This means that the
pull-back surface has more than one non-reduced fibre. Then there
is a quadratic twist of the surface, sending $x\mapsto \alpha^2 x$
and $y\mapsto \alpha^3y$, which replaces an even number of
non-reduced fibres by their reduced relatives (i.e.~$I_n^{\ast}$
by $I_n$ and $II^{\ast},III^{\ast},IV^{\ast}$ by $IV,III,II$,
respectively). Here, $\alpha$ is the vanishing polynomial of the
cusps of the non-reduced fibres. Following \cite{M2} this process
will be called \emph{deflation}. Since our main interest lies in
elliptic (K3) surfaces with only semi-stable fibres, deflation
provides a useful tool to construct such surfaces via base change.

It is exactly these two methods (minimalization and deflation after a suitable base change) which we will use to resolve the non-semi-stable fibres of %%@
the base surface $Y$. The explicit behaviour of the singular fibres under a base change can be derived from \cite[IV.9, table 4.1]{Si} or found in %%@
\cite[section 7]{MP1}. We sketch it in the next figure where the number next to an arrow %%@
denotes the order of ramification under $\pi$ (of one particular
pre-image). The fibres of type $I_n^{\ast}$ are exceptional in
that they allow two possibilities of substitution by semi-stable
fibres: Either by ramification of even index or by pairwise
deflation.

\begin{figure}[!ht]
\begin{center}
\setlength{\unitlength}{0.08mm}
\begin{picture}(1730,370)

\put(200,300){\circle*{10}}
\put(500,300){\circle*{10}}
\put(800,300){\circle*{10}}
\put(1100,300){\circle*{10}}

\put(200,100){\circle*{10}}
\put(500,100){\circle*{10}}
\put(800,100){\circle*{10}}
\put(1100,100){\circle*{10}}

\put(201,300){\line(0,-1){200}}
\put(500,300){\line(0,-1){200}}
\put(800,300){\line(0,-1){200}}
\put(1100,300){\line(0,-1){200}}

\put(193,320){$I_{2n}$}
\put(192,50){$I_n^{\ast}$}
\put(493,320){$I_0^{\ast}$}
\put(488,50){$II^{(\ast)}$}
\put(793,320){$I_0^{\ast}$}
\put(782,50){$III^{(\ast)}$}
\put(1095,320){$I_0$}
\put(1084,50){$IV^{(\ast)}$}

\put(215,195){2}
\put(515,195){3}
\put(815,195){2}
\put(1115,195){3}

\put(1440,300){\vector(0,-1){200}}
\put(1455,195){$\pi$}
\put(1428,320){$X$}
\put(1428,50){$Y$}

\end{picture}
\end{center}
\caption{The resolution of the non-semi-stable fibres}
\label{Fig:res}
\end{figure}
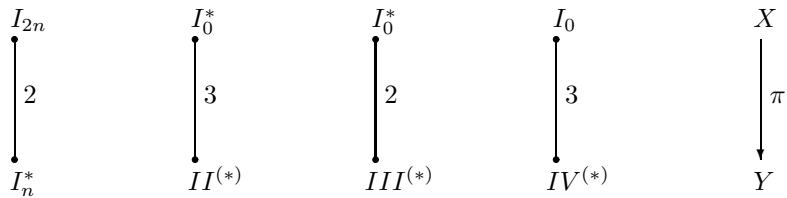

\section{The deflated base changes}
\label{s:deflated}

Our main interest lies in finding equations over $\Q$ for
\emph{extremal} semi-stable elliptic K3 fibrations. By definition,
these are singular elliptic K3 surfaces with finite Mordell-Weil
group and only semi-stable fibres. These assumptions are quite
restrictive. It is an immediate consequence that the number of
cusps has to be 6, and one finds the 112 possible configurations
of singular fibres in the classification of \cite[Thm (3.1)]{MP2}.
The main idea of this paper is to produce some of these K3
fibrations by the methods described in the previous section via
the pull-back of a rational elliptic surface by a base change.
This approach is greatly helped by the good explicit knowledge one
has of the rational elliptic surfaces (cf.~\cite{H}, \cite{MP1},
\cite{S-H}). One only has to construct suitable base changes.

The starting point for our considerations is a rational elliptic
surface with a section
\[
Y \stackrel{r}{\rightarrow} \Phoch{1}.
\]
For a base change $\pi:
\Phoch{1}\rightarrow\Phoch{1}$, let
\[
X \stackrel{\pi\circ r}{\longrightarrow}\Phoch{1}
\]
denote the pullback via $\pi$. Since the Mordell-Weil group
$MW(Y)$ injects into $MW(X)$, we will assume $Y$ to be extremal.
Note, however, that the process of deflation can a priori change
the Mordell-Weil group. Therefore it could also seem worth
considering non-extremal rational elliptic surfaces, especially
those with a small number of cusps as presented in \cite{H} and
\cite{S-H}. A close observation nevertheless shows that these
would not produce any configurations different from those known or
obtained in this paper, unless one turns to the general case where
$\pi$ has degree 24. This can also be derived from
\cite[Thm.~1.2]{Kl}. A solution to the general case has recently
been announced by Beukers and Montanus \cite{BM}.

Extremal rational elliptic surfaces have been completely
classified by Miranda-Persson in \cite{MP1}. There are six
semi-stable surfaces with four cusps. These had previously been
identified as Shioda modular by Beauville \cite{B}. As explained,
these surfaces have been treated exhaustively in \cite{TY}, giving
rise to the semi-stable extremal elementary fibrations of
\cite{P}.

Furthermore there are four surfaces with only two cusps (one of
them appearing in a continuous family). We will not use these
surfaces since they have no fibre of type $I_n$ or $I_n^{\ast}$
with $n>0$ at all.

The remaining six extremal rational elliptic surfaces have three
cusps. Each has exactly one non-reduced fibre while the other two
singular fibres are semi-stable. These are the surfaces we are
going to investigate for a pull-back via a base change. For the
remaining part of this section we will concentrate on those
("deflated") base changes $\pi$ which give rise to a non-inflated
pull-back K3 surface $X \stackrel{\pi\circ
r}{\longrightarrow}\Phoch{1}$ after minimalizing. The pull-back
surfaces coming from inflating base changes will be derived in
Section \ref{s:inflating}.

Let $Y$ be an extremal rational elliptic surface $Y$ with three
cusps. There are a number of conditions on the base change $\pi$
if one wants an extremal semi-stable K3 pull-back. The Euler
number $e(X)=24$ of the pull-back surface $X$ predicts the degree
of $\pi$, only depending on the type of the non-reduced fibre
$W^{\ast}$ of $Y$.

On the one hand, if $W$ is of additive type
(i.e.~$W\in\{II,III,IV\}$), we will eventually replace it by
smooth fibres after minimalizing and deflating, if necessary. Let
the other two singular fibres be of type $I_m$ and $I_n$ with
$m,n\in\N$. Then $m+n\leq 4$, and $e(X)=$ (deg $\pi)(m+n)$. On the
other hand, let the non-reduced fibre have type $I_k^{\ast}$. If
the other singular fibres have again $m$ and $n$ components, then
we have $k+m+n=6$. Therefore we require deg $\pi=4$.

Our assumption, that the semi-stable pull-back $X$ is extremal,
implies the minimal number of six singular fibres. This gives
another stringent restriction. In some cases, this leads to a
contradiction to the Hurwitz formula
\[
-2\geq -2~\text{deg}~\pi + \sum_{x\in\Phoch{1}}(\text{deg}~\pi - \#\pi^{-1}(x)).
\]
Finally, as we decided to concentrate on resolving the non-reduced
fibre $W^{\ast}$ by a deflated base change, the ramification index
at the corresponding cusp has to be divisible by 2, 3, 4, or 6 if
$W=I_n^*, IV, III,$ or $II$, respectively (cf.~Fig.
\ref{Fig:res}).

It will turn out that some of the base changes can only be defined
over an extension of $\Q$ of low degree. Nevertheless, for any
base change, the pull-back surface $X$ will have at least two
rational cusps. For simplicity and without loss of generality, we
will choose these by a M\"obius transformation to be 0 and
$\infty$ (and a further third rational cusp, if it exists, to be
1). For every rational surface, we will derive a Weierstrass
equation over $\Q$ such that the cusps are 0,1 and $\infty$. This
gives us the opportunity to construct the base changes before
considering the surfaces.

We require that the pull-back of a base change $\pi$ does not
factor through a Shioda modular rational elliptic surface.
Equivalently, $\pi$ does not factor into a composition
$\pi''\circ\pi'$ of a degree 2 map $\pi''$  and a further map
$\pi'$, such that the non-reduced fibre is already resolved by
$\pi'$. Otherwise, the intermediate pull-back $X'
\stackrel{\pi'\circ r}{\longrightarrow}\Phoch{1}$ would be
semi-stable and have at most 4 cusps, hence it would be Shioda
modular by \cite{B}.

We shall now investigate the deflated base changes with the above
listed properties. For an extremal rational elliptic surface $Y$
with three cusps, our analysis depends only on the type of the
non-reduced fibre $W^{\ast}$. Throughout we employ the notation of
\cite{MP1}. For the computations we wrote a straight forward Maple
program. In all but two cases, this sufficed to determine the
minimal base change. For the remaining two base changes, we used
Macaulay to compute a solution mod $p$ for some small primes $p$
and then lift to characteristic 0.

At first assume $\mathbf{W}$ to be $\mathbf{II}$. According to
\cite[section 5]{MP1} there is, up to isomorphism, a unique
rational elliptic surface with this fibre whose other two singular
fibres are both of type $I_1$. Denote this by $\mathbf{X_{211}}$.
For the pull-back $X$ to have Euler number $e(X)=24$, we would
need $\pi$ to have degree 12 and ramification index 12 or (6,6) at
the cusp of the non-reduced fibre. Then, the restriction of the
other two cusps to have exactly six pre-images under $\pi$ leads
to a contradiction to the Hurwitz formula.

The situation is similar if $\mathbf{W=III}$. This implies
$Y=\mathbf{\mathbf{X_{321}}}$ to have further singular fibres
$I_2, I_1$. The $III^{\ast}$ fibre requires ramification of index
a multiple of 4, so $\pi$ must have degree 8 with ramification
index 8 or (4,4) at the cusp of this fibre. Again, the Hurwitz
formula rules out the pull-back surface $X$ to have only six
(semi-stable) singular fibres. In Section \ref{s:inflating}, we
will construct inflating base changes for this surface which
resolve the non-reduced fibre and are compatible with the Hurwitz
formula.

$\mathbf{W=IV}$ gives a priori two possibilities for the elliptic
surface $Y$. One of them, $\mathbf{X_{431}}$ with singular fibres
of type $IV^{\ast}, I_3, I_1$, does actually exist by \cite{MP1}.
Since a $IV^{\ast}$ fibre requires the ramification index to be
divisible by 3, we need $\pi$ to have degree 6 and ramification of
order (3,3) at the cusp of this fibre (since ramification index 6
would contradict the other two fibres having six pre-images again
by the Hurwitz formula). The suitable maps are presented in the
next paragraphs. Throughout we assume the $IV^{\ast}$ fibre to sit
above 1. If such a base change was totally ramified above one of
the two remaining cusps, then it would necessarily be composite.
Since this was excluded, we only have to deal with those maps such
that 0 has two or three pre-images (and then exchange 0 and
$\infty$).

Consider a base change $\pi$ of degree 6 with ramification index
(3,3) at 1. At first, let 0 and $\infty$ have three pre-images.
Our restriction that $\pi$ does not factor, implies that at least
one of the cusps has ramification index (3,2,1). We will assume
$\infty$ to have this ramification index. Then we search for maps
such that 0 has ramification index (4,1,1), (3,2,1) or (2,2,2).
However, a map with the last ramification cannot exist since the
resulting pull-back of $\mathbf{X_{431}}$ does not appear in the
list of \cite{MP2}. In the following, we will frequently use this
argument. Here, let us once give the details.

Assume that such a base change $\pi$ exists. Then we can realize
the configuration [2,2,2,3,6,9] as pull-back $X$ from
$\mathbf{X_{431}}$. We claim that this fibration has a 2-section.
Otherwise, the fibre types would imply that
\[
(\Z/2)^4\subseteq NS(X)^\vee/NS(X).
\]
So the 2-length of the quotient would be at least 4. Let $T_X$
denote the transcendental lattice of $X$, that is
$T_X=NS(X)^\bot\subset\text{H}^2(X,\Z)$. Then by lattice theory
\cite[\S 1]{N},
\[
NS(X)^\vee/NS(X)\cong T_X^\vee/T_X.
\]
Since here $T_X$ has rank 2, this quotient can maximally have 2-length 2. This
gives the required contradiction.

So $X$ has a 2-section $\sigma$. Consider the quotient of $X$ by
(translation by) $\sigma$. The minimal desingularization $Z$ of
this is again an elliptic K3 surface. The resulting singular
fibres can be computed from the components which $\sigma$ meets.
In particular, the fibres of types $I_9$ and $I_3$ result in
$I_{18}$ and $I_6$ in $Z$. But this is impossible, since then
$e(Z)>24$. Hence, $X$ cannot admit a 2-section. Since the
existence of the 2-section followed from the construction via
$\pi$, this base change cannot exist.

The computations show that the base change with the second
ramification, $\mathbf{\tilde\pi}$, can only be realized over a
cubic extension of $\Q$. (With $v$ a solution of
$5x^3+12x^2+12x+4$ it can be given as
$\mathbf{\tilde{\pi}}((s:t))=(s^3(s-t)^2(s+(2+3v)t):-(2+3v)t^3(s+(1+v)^2t)^2(s+t/(5v+2)))$.
Here, we only construct the first base change:
\begin{eqnarray*}
\mathbf{\pi_{3,4}}:~~~ \Phoch{1}~ & \rightarrow & \quad\Phoch{1}\\
(s:t) & \mapsto & (27s^4(125t^2-90st-27s^2):-3125t^3(t-s)^2(5t+4s)).
\end{eqnarray*}
Then
$27s^4(125t^2-90st-27s^2)+3125t^3(t-s)^2(5t+4s)=(25t^2-10st-9s^2)^3$,
so $\pi_{3,4}$ has the required properties.

Now, let 0 have only two pre-images and $\infty$ four. The
respective ramification indices are (5,1), (4,2) or (3,3) and
(2,2,1,1) or (3,1,1,1). Only the first and the last of these do
not allow a composition as a degree 3 and a degree 2 map. Hence,
at least one of these two ramification indices must occur for the
base change to meet our criteria.

Let us first construct the maps with ramification index (5,1) at
0. They are:
\begin{eqnarray*}
\mathbf{\pi_{5,3}}:~~~ \Phoch{1}~ & \rightarrow & \quad\Phoch{1}\\
(s:t) & \mapsto & (729s^5(s-t):-t^3(135s^3+9st^2+t^3))
\end{eqnarray*}
with $729s^5(s-t)+t^3(135s^3+9st^2+t^3)=(9s^2-3st-t^2)^3$ and
\begin{eqnarray*}
\mathbf{\pi_{5,2}}:~~~ \Phoch{1}~ & \rightarrow & \quad\Phoch{1}\\
(s:t) & \mapsto & (2^63^3s^5t:-(s^2-4st-t^2)^2(125s^2+22st+t^2))
\end{eqnarray*}
with $2^63^3s^5t+(s^2-4st-t^2)^2(125s^2+22st+t^2)=(5s^2+10st+t^2)^3.$

The other base changes have ramification index (3,1,1,1) at
$\infty$. It is immediate that there is no such map $\pi$ with
ramification index (3,3) at 0: After exchanging 1 and $\infty$,
the map $\pi$ would have to look like $(f_0^3g_0^3:f_1^3g_1^3)$
with distinct linear homogeneous factors $f_i,g_i$. Then, with
$\varrho$ a primitive third root of unity,
\[f_0^3g_0^3-f_1^3g_1^3=(f_0g_0-f_1g_1)(f_0g_0-\varrho
f_1g_1)(f_0g_0-\varrho^2 f_1g_1).\] This polynomial cannot have a
cubic factor. Therefore, the next map completes the list of
suitable base changes for $\mathbf{X_{431}}$. Remember that we
still have to take the permutation of 0 and $\infty$ via
exchanging $s$ and $t$ into account.
\begin{eqnarray*}
\mathbf{\pi_{4,3}}:~~~ \Phoch{1}~ & \rightarrow & \quad\Phoch{1}\\
(s:t) & \mapsto & (729s^4t^2:-(s-t)^3(8s^3+120s^2t-21st^2+t^3))
\end{eqnarray*}
with $729s^4t^2+(s-t)^3(8s^3+120s^2t-21st^2+t^3)=(2 s^2-8 s t-t^2)^3$.

\vspace{0.4cm}

We conclude this section by considering the non-reduced fibre
$\mathbf{W^{\ast}}$ to equal $\mathbf{I_n^{\ast}}$ for some $n>
0$. By \cite{MP1} there are three extremal rational elliptic
surfaces with such a singular fibre. All of them have two further
singular fibres, both semi-stable. The surfaces will be introduced
in the next section. Independent of the surface, we have already
seen that an adequate deflated base change $\pi$ must have degree
4 and ramification of index (2,2) or 4 at the cusp of the
$I_n^{\ast}$ fibre. Then $X$ is extremal if and only if the two
other cusps have 4 or 5 pre-images, respectively. Assume that the
non-reduced fibre sits above the cusp $\infty=(1:0)$. We now
construct the base changes $\pi$ which do not factor into two maps
of degree 2. Equivalently, one of the cusps has ramification index
(3,1).

At first, let the base change $\pi$ be totally ramified at
$\infty$. By the above considerations, the other two cusps have
ramification indices (3,1) and (2,1,1). Up to exchanging them, for
example by
\[
\phi: (s:t) \mapsto (t-s:t),
\]
the map $\pi$ can be realized as
\begin{eqnarray*}
\mathbf{\pi_4}:~~ \Phoch{1} ~& \rightarrow & \quad\Phoch{1}\\
(s:t) & \mapsto & (256s^3(s-t):-27t^4),
\end{eqnarray*}
since %%@
$(256s^3(s-t)+27t^4)=(4s-3t)^2(16s^2+8st+3t^2)$.

In the other case, $\pi$ has ramification index (2,2) at $\infty$.
Our restrictions imply ramification of index (3,1) at one of the
other two cusps. Without loss of generality, let this cusp be 1.
Then the last cusp also has two pre-images and thus ramification
index (3,1) or (2,2). The second cannot exist: Given such a map,
it could be expressed as $(f_0^2g_0^2:f_1^2g_1^2)$ with distinct
homogeneous linear forms $f_i,g_i$ in $s,t$. The factorization
$f_0^2g_0^2-f_1^2g_1^2=(f_0g_0+f_1g_1)(f_0g_0-f_1g_1)$ cannot have
a cubic factor. Hence 1 cannot have ramification index (3,1).

We conclude this section with a base change $\mathbf{\pi_2}$ of
degree 4 with ramification indices (3,1), (3,1) and (2,2):
\begin{eqnarray*}
\pi_2:~~ \Phoch{1} ~& \rightarrow & \quad\Phoch{1}\\
(s:t) & \mapsto & (64s^3(s-t):(8s^2-4st-t^2)^2).
\end{eqnarray*}

In the next two sections (\ref{s:degree 4}, \ref{s:degree 6}), we
will substitute the base changes $\pi_{3,4}, \pi_{5,3}, \pi_{5,2}, \pi_{4,3}, \pi_2$ and $\pi_4$ into the normalized
Weierstrass equations of the extremal rational elliptic surfaces
with three singular fibres. We will derive equations over $\Q$ for
extremal K3 surfaces with six singular fibres, all of which
semi-stable.

\section{The fibrations coming from degree 4 base changes}
\label{s:degree 4}

In this and the next section (\ref{s:degree 6}), we proceed as
follows to obtain equations for extremal K3 surface with six
semi-stable fibres: Consider the Weierstrass equations given for
the extremal rational elliptic surfaces in \cite[Table 5.2]{MP1}.
We apply the normalizing M\"obius transformation which maps the
cusps to 0, 1 and $\infty$. Then we exhibit the deflated base
changes $\pi_i$ from the previous section. After minimalizing by
an admissible change of variables, this gives the Weierstrass
equations for 16 extremal elliptic K3 surfaces from the list of
\cite{MP2}. Throughout we choose the coefficients of the
polynomials $A$ and $B$ involved in the Weierstrass equation to be
minimal by rescaling, if necessary. By construction, the pull-back
surface $X$ always inherits the sections of the rational elliptic
surface $Y$. As a consequence, we are able to derive the
isomorphism class of $X$ (in terms of the intersection form on its
transcendental lattice) from the classification in \cite{SZ}.

In this section we consider only the extremal rational elliptic
surfaces with an $\mathbf{I_n^{\ast}}$ fibre ($n>0$). As
explained, they require a base change of degree 4. Before
substituting by the base changes $\pi_4$ or $\pi_2$, we choose the
normalizing M\"obius transformation in such a way that the
$I_n^{\ast}$ fibre sits above $\infty$.

Let us start with $\mathbf{X_{4 1 1}}$ which has Weierstrass equation
\[
y^2 = x^3-3\, t^2 (s^2-3 t^2)\, x+s t^3 (2 s^2-9 t^2)
\]
The singular fibres are $I_4^{\ast}$ over $\infty$ and two $I_1$
over $\pm 2$. Substituting $(s:t)\mapsto(4 s-2t:t)$ maps the two
$I_1$ fibres to 0 and 1. This gives
\[
y^2 = x^3-3\, t^2 (16 s^2-16 s t+t^2)\,x+2\, t^3 (2 s-t) (32 s^2-32 s t-t^2).
\]

Then, we substitute by $\mathbf{\pi_4}$  and get the Weierstrass equation %%@
of an extremal K3 surface:
\begin{small}
\begin{eqnarray*}
y^2 & = & x^3-3\,(9s^8+48s^7t+48s^4t^4+64s^6t^2+128s^3t^5+16t^8)\,x\\
&& \quad-2\,(3s^4+8s^3t+8t^4)(9s^8+48s^7t+48s^4t^4+64s^6t^2+128s^3t^5-8t^8).
\end{eqnarray*}
\end{small}
This provides a realization of the configuration
\textbf{[1,1,1,2,3,16]} in the notation of \cite{MP2} (i.e.~3
fibres of type $I_1$, one of types $I_2, I_3$ and $I_{16}$ each).

%\pi_4:
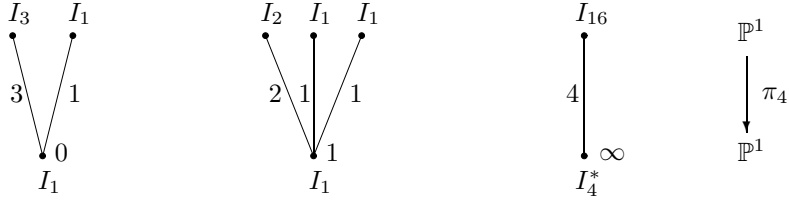
\begin{figure}[!ht]
\begin{center}
\setlength{\unitlength}{0.08mm}
\begin{picture}(1500,320)

\put(180,265){\circle*{10}}
\put(172,290){$I_1$}
\put(180,265){\line(-1,-4){50}}
\put(172,155){1}
\put(130,65){\circle*{10}}
\put(150,56){0}
\put(122,10){$I_1$}
\put(80,265){\circle*{10}}
\put(72,290){$I_3$}
\put(80,265){\line(1,-4){50}}
\put(76,155){3}

\put(580,65){\circle*{10}}
\put(572,10){$I_1$}
\put(600,56){1}
\put(580,265){\circle*{10}}
\put(572,290){$I_1$}
\put(581,265){\line(0,-1){200}}
\put(555,155){1}
\put(500,265){\circle*{10}}
\put(492,290){$I_2$}
\put(500,265){\line(2,-5){80}}
\put(506,155){2}
\put(660,265){\circle*{10}}
\put(652,290){$I_1$}
\put(660,265){\line(-2,-5){80}}
\put(640,155){1}

\put(1030,65){\circle*{10}}
\put(1015,10){$I_4^{\ast}$}
\put(1030,265){\circle*{10}}
\put(1030,265){\line(0,-1){200}}
\put(1000,155){4}
\put(1055,58){$\infty$}
\put(1015,290){$I_{16}$}

\put(1285,257){$\Phoch{1}$}
\put(1300,230){\vector(0,-1){125}}
\put(1285,57){$\Phoch{1}$}
\put(1325,160){$\pi_4$}

\end{picture}
\end{center}
\caption{A realization of \textbf{[1,1,1,2,3,16]}}
\label{fig:24}
\end{figure}
\vspace{0.2cm}

On the other hand, we can also substitute by $\mathbf{\pi_2}$ in
the normalized Weierstrass equation. We obtain:
\begin{footnotesize}
\begin{eqnarray*}
y^2 & = & x^3-3\, (16 s^8-64 s^7 t-224 s^5 t^3+392 t^4 s^4+64 s^6 t^2+112 t^5 s^3+16 t^6 s^2+8 t^7 s+t^8)\,x\\
&&\quad-2\, (2 s^2-4 s t-t^2) (2 s^2+t^2)\\
&&\quad\quad (16 s^8-64 s^7 t+544 s^5 t^3-952 t^4 s^4+64 s^6
t^2-272 t^5 s^3+16 t^6 s^2+8 t^7 s+t^8).
\end{eqnarray*}
\end{footnotesize}
This realizes \textbf{[1,1,3,3,8,8]}.

%\pi_2:
\begin{figure}[!ht]
\begin{center}
\setlength{\unitlength}{0.08mm}
\begin{picture}(1500,320)

\put(180,265){\circle*{10}}
\put(172,290){$I_1$}
\put(180,265){\line(-1,-4){50}}
\put(172,155){1}
\put(130,65){\circle*{10}}
\put(150,56){0}
\put(122,10){$I_1$}
\put(80,265){\circle*{10}}
\put(72,290){$I_3$}
\put(80,265){\line(1,-4){50}}
\put(76,155){3}

\put(610,265){\circle*{10}}
\put(602,290){$I_1$}
\put(610,265){\line(-1,-4){50}}
\put(602,155){1}
\put(560,65){\circle*{10}}
\put(580,56){1}
\put(552,10){$I_1$}
\put(510,265){\circle*{10}}
\put(502,290){$I_3$}
\put(510,265){\line(1,-4){50}}
\put(506,155){3}

\put(1040,265){\circle*{10}}
\put(1032,290){$I_8$}
\put(1040,265){\line(-1,-4){50}}
\put(1030,155){2}
\put(990,65){\circle*{10}}
\put(1015,58){$\infty$}
\put(975,10){$I_4^{\ast}$}
\put(940,265){\circle*{10}}
\put(932,290){$I_8$}
\put(940,265){\line(1,-4){50}}
\put(936,155){2}

\put(1285,257){$\Phoch{1}$}
\put(1300,230){\vector(0,-1){125}}
\put(1285,57){$\Phoch{1}$}
\put(1325,160){$\pi_2$}

\end{picture}
\end{center}
\caption{A realization of \textbf{[1,1,3,3,8,8]}}
\end{figure}
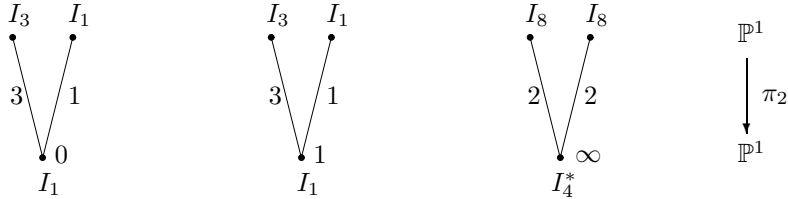
\vspace{0.2cm}

We shall apply the same procedure to the surface $\mathbf{X_{1 4
1}}$ with Weierstrass equation
\[
y^2 = x^3-3\,(s-2 t)^2 (s^2-3 t^2)\, x+s (s-2 t)^3 (2 s^2-9 t^2).
\]
The normalization of the cusps leads to the
Weierstrass equation
\[
y^2=x^3-3\, t^2 (16 t^2-16 s t+s^2)\,x+2 \,t^3 (s-2 t)
(s^2+32st-32 t^2).
\]
This has an $I_4$ fibre above 0, $I_1$ above 1 and $I_1^*$ above
$\infty$. From this, we will derive three examples, two of them
connected by the M\"obius transformation $\phi$. Substitution of
$\mathbf{\pi_4}$ produces a realization of
\textbf{[1,1,2,4,4,12]}:
\begin{small}
\begin{eqnarray*}
y^2 & = & x^3-3\,(16t^8+48s^4t^4-64s^3t^5+9s^8-24s^7t+16s^6t^2)\,x\\
&& \quad-2\,(2t^4+3s^4-4s^3t)(-32t^8-96s^4t^4+128s^3t^5+9s^8-24s^7t+16s^6t^2)
\end{eqnarray*}
\end{small}

%\pi_4:
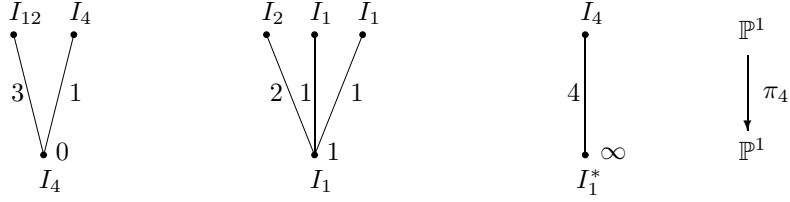
\begin{figure}[!ht]
\begin{center}
\setlength{\unitlength}{0.08mm}
\begin{picture}(1500,320)

\put(180,265){\circle*{10}}
\put(172,290){$I_4$}
\put(180,265){\line(-1,-4){50}}
\put(172,155){1}
\put(130,65){\circle*{10}}
\put(150,56){0}
\put(122,10){$I_4$}
\put(80,265){\circle*{10}}
\put(72,290){$I_{12}$}
\put(80,265){\line(1,-4){50}}
\put(76,155){3}

\put(580,65){\circle*{10}}
\put(572,10){$I_1$}
\put(600,56){1}
\put(580,265){\circle*{10}}
\put(572,290){$I_1$}
\put(581,265){\line(0,-1){200}}
\put(555,155){1}
\put(500,265){\circle*{10}}
\put(492,290){$I_2$}
\put(500,265){\line(2,-5){80}}
\put(506,155){2}
\put(660,265){\circle*{10}}
\put(652,290){$I_1$}
\put(660,265){\line(-2,-5){80}}
\put(640,155){1}

\put(1030,65){\circle*{10}}
\put(1015,10){$I_1^{\ast}$}
\put(1030,265){\circle*{10}}
\put(1030,265){\line(0,-1){200}}
\put(1000,155){4}
\put(1055,58){$\infty$}
\put(1022,290){$I_4$}

\put(1285,257){$\Phoch{1}$}
\put(1300,230){\vector(0,-1){125}}
\put(1285,57){$\Phoch{1}$}
\put(1325,160){$\pi_4$}

\end{picture}
\end{center}
\caption{A realization of \textbf{[1,1,2,4,4,12]}}
\end{figure}
\vspace{0.2cm}

Conjugation by $\phi$ gives the Weierstrass equation
\[
y^2=x^3-3\, t^2 (t^2+14 s t+s^2)\,x-2\, (t+s) t^3 (t^2-34 s
t+s^2).
\]
This leads to the following realization of \textbf{[1,3,4,4,4,8]}
\begin{small}
\begin{eqnarray*}
y^2 & = & x^3+3\,(-24s^7t-16s^6t^2-9s^8-t^8+42s^4t^4+56s^3t^5)\,x\\
&& \quad -2\,(9s^8+24s^7t+16s^6t^2+102s^4t^4+136s^3t^5+t^8)(-3s^4-4s^3t+t^4)
\end{eqnarray*}
\end{small}

%\pi_4:
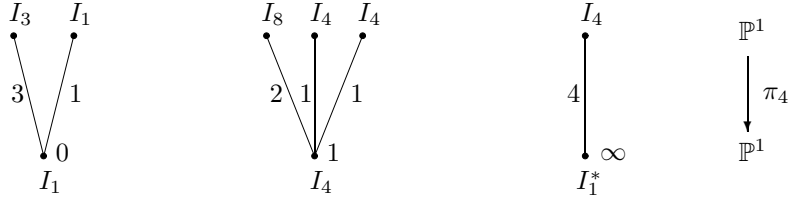
\begin{figure}[!ht]
\begin{center}
\setlength{\unitlength}{0.08mm}
\begin{picture}(1500,320)

\put(180,265){\circle*{10}}
\put(172,290){$I_1$}
\put(180,265){\line(-1,-4){50}}
\put(172,155){1}
\put(130,65){\circle*{10}}
\put(150,56){0}
\put(122,10){$I_1$}
\put(80,265){\circle*{10}}
\put(72,290){$I_3$}
\put(80,265){\line(1,-4){50}}
\put(76,155){3}

\put(580,65){\circle*{10}}
\put(572,10){$I_4$}
\put(600,56){1}
\put(580,265){\circle*{10}}
\put(572,290){$I_4$}
\put(581,265){\line(0,-1){200}}
\put(555,155){1}
\put(500,265){\circle*{10}}
\put(492,290){$I_8$}
\put(500,265){\line(2,-5){80}}
\put(506,155){2}
\put(660,265){\circle*{10}}
\put(652,290){$I_4$}
\put(660,265){\line(-2,-5){80}}
\put(640,155){1}

\put(1030,65){\circle*{10}}
\put(1015,10){$I_1^{\ast}$}
\put(1030,265){\circle*{10}}
\put(1030,265){\line(0,-1){200}}
\put(1000,155){4}
\put(1055,58){$\infty$}
\put(1022,290){$I_4$}

\put(1285,257){$\Phoch{1}$}
\put(1300,230){\vector(0,-1){125}}
\put(1285,57){$\Phoch{1}$}
\put(1325,160){$\pi_4$}

\end{picture}
\end{center}
\caption{A realization of \textbf{[1,3,4,4,4,8]}}
\end{figure}
\vspace{0.2cm}

On the other hand, substitution of $\mathbf{\pi_2}$ in the normalized Weierstrass equation gives:
\begin{small}
\begin{eqnarray*}
y^2 & = & x^3-3\,(s^8-4 s^7 t+16 s^5 t^3-28 t^4 s^4+4 s^6 t^2-8 t^5 s^3+16 s^2 t^6+8 s t^7+t^8)\,x
\\
&& \quad -(2 s^4-4 s^3 t+4 s t^3+t^4)\\
&& \quad\quad (s^8-4 s^7 t-32 s^5 t^3+56 t^4 s^4+4 s^6 t^2+16 t^5 s^3-32 s^2 t^6-16 s t^7-2 t^8).
\end{eqnarray*}
\end{small}
This realization of \textbf{[1,2,2,3,4,12]} is invariant under
conjugation by $\phi$.

%\pi_2:
\begin{figure}[!ht]
\begin{center}
\setlength{\unitlength}{0.08mm}
\begin{picture}(1500,320)

\put(180,265){\circle*{10}}
\put(172,290){$I_4$}
\put(180,265){\line(-1,-4){50}}
\put(172,155){1}
\put(130,65){\circle*{10}}
\put(150,56){0}
\put(122,10){$I_4$}
\put(80,265){\circle*{10}}
\put(65,290){$I_{12}$}
\put(80,265){\line(1,-4){50}}
\put(76,155){3}

\put(610,265){\circle*{10}}
\put(602,290){$I_1$}
\put(610,265){\line(-1,-4){50}}
\put(602,155){1}
\put(560,65){\circle*{10}}
\put(580,56){1}
\put(552,10){$I_1$}
\put(510,265){\circle*{10}}
\put(502,290){$I_3$}
\put(510,265){\line(1,-4){50}}
\put(506,155){3}

\put(1040,265){\circle*{10}}
\put(1032,290){$I_2$}
\put(1040,265){\line(-1,-4){50}}
\put(1030,155){2}
\put(990,65){\circle*{10}}
\put(1015,58){$\infty$}
\put(975,10){$I_1^{\ast}$}
\put(940,265){\circle*{10}}
\put(932,290){$I_2$}
\put(940,265){\line(1,-4){50}}
\put(936,155){2}

\put(1285,257){$\Phoch{1}$}
\put(1300,230){\vector(0,-1){125}}
\put(1285,57){$\Phoch{1}$}
\put(1325,160){$\pi_2$}

\end{picture}
\end{center}
\caption{A realization of \textbf{[1,2,2,3,4,12]}}
\label{Fig:1,2,2,3,4,12}
\end{figure}
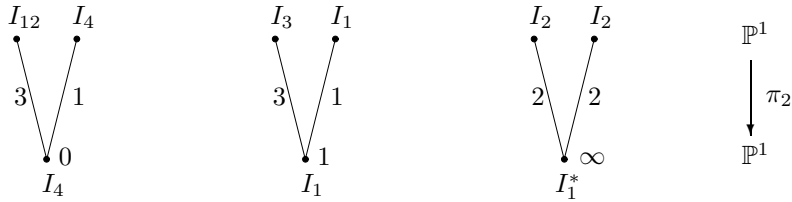
\vspace{0.2cm}

Finally, we turn to the surface $\mathbf{X_{222}}$ in the notation
of \cite{MP1}. Miranda-Persson give the Weierstrass equation
\[
y^2=x^3-3\, s t (s - t)^2\,x+(s-t)^3(s^3+t^3).
\]
This has cusps the third roots of unity with an $I_2^{\ast}$ fibre
above 1 and $I_2$ fibres above the two primitive roots $\omega,
\omega^2$. Take $\omega$ in the upper half plane. Consider the
M\"obius transformation which maps $\omega$ to $\infty$ and
$\omega^2$ to 0 while fixing 1. This gives a Weierstrass equation
which is not defined over $\Q$. However, the change of variables
$x\mapsto \xi^2x, y\mapsto \xi^3y$ with $\xi=3\sqrt{-3}$ leads to
a Weierstrass equation over $\Q$ with the same cusps:
\[
y^2=x^3-3\,(s^2-st+t^2)(s-t)^2x+(s-2t)(2s-t)(t+s)(s-t)^3.
\]

We exchange the cusps 1 and $\infty$ and subsequently substitute
by $\mathbf{\pi_4}$ or $\mathbf{\pi_2}$. The first substitution
gives a realization of \textbf{[2,2,2,4,6,8]}:
\begin{eqnarray*}
y^2 & = & x^3-3\,(9s^8-24s^7t+16s^6t^2+3s^4t^4-4s^3t^5+t^8)\,x\\
&& \quad-(2t^4+3s^4-4s^3t)(3s^4-4s^3t-t^4)(6s^4-8s^3t+t^4)
\end{eqnarray*}

%\pi_4:
\begin{figure}[!ht]
\begin{center}
\setlength{\unitlength}{0.08mm}
\begin{picture}(1500,320)

\put(180,265){\circle*{10}}
\put(172,290){$I_2$}
\put(180,265){\line(-1,-4){50}}
\put(172,155){1}
\put(130,65){\circle*{10}}
\put(150,56){0}
\put(122,10){$I_2$}
\put(80,265){\circle*{10}}
\put(72,290){$I_6$}
\put(80,265){\line(1,-4){50}}
\put(76,155){3}

\put(580,65){\circle*{10}}
\put(572,10){$I_2$}
\put(600,56){1}
\put(580,265){\circle*{10}}
\put(572,290){$I_2$}
\put(581,265){\line(0,-1){200}}
\put(555,155){1}
\put(500,265){\circle*{10}}
\put(492,290){$I_4$}
\put(500,265){\line(2,-5){80}}
\put(506,155){2}
\put(660,265){\circle*{10}}
\put(652,290){$I_2$}
\put(660,265){\line(-2,-5){80}}
\put(640,155){1}

\put(1030,65){\circle*{10}}
\put(1015,10){$I_2^{\ast}$}
\put(1030,265){\circle*{10}}
\put(1030,265){\line(0,-1){200}}
\put(1000,155){4}
\put(1055,58){$\infty$}
\put(1022,290){$I_8$}

\put(1285,257){$\Phoch{1}$}
\put(1300,230){\vector(0,-1){125}}
\put(1285,57){$\Phoch{1}$}
\put(1325,160){$\pi_4$}

\end{picture}
\end{center}
\caption{A realization of \textbf{[2,2,2,4,6,8]}}
\label{Fig:2,2,2,4,6,8}
\end{figure}
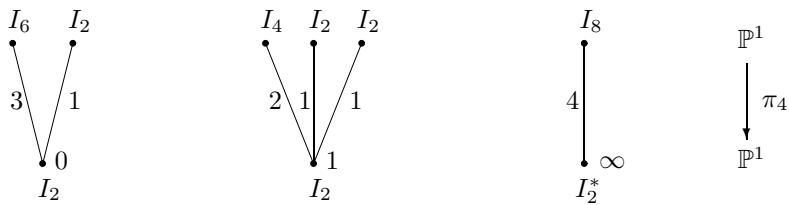
\vspace{0.2cm}

The second substitution realizes \textbf{[2,2,4,4,6,6]}:
\begin{small}
\begin{eqnarray*}
y^2 & = & x^3-3\,(16 s^8-64 t s^7+64 t^2 s^6+16 t^3 s^5-28 t^4 s^4-8 t^5 s^3+16 t^6 s^2+8 t^7 s+t^8)\,x\\
&& \quad +2\,(2s^2-4st-t^2)(8 s^4-16 s^3 t+4 s t^3+t^4)(t^2+2 s^2)(2 s^4-4 s^3 t+4 s t^3+t^4) .
\end{eqnarray*}
\end{small}

%\pi_2:
\begin{figure}[!ht]
\begin{center}
\setlength{\unitlength}{0.08mm}
\begin{picture}(1500,320)

\put(180,265){\circle*{10}}
\put(172,290){$I_2$}
\put(180,265){\line(-1,-4){50}}
\put(172,155){1}
\put(130,65){\circle*{10}}
\put(150,56){0}
\put(122,10){$I_2$}
\put(80,265){\circle*{10}}
\put(72,290){$I_6$}
\put(80,265){\line(1,-4){50}}
\put(76,155){3}

\put(610,265){\circle*{10}}
\put(602,290){$I_2$}
\put(610,265){\line(-1,-4){50}}
\put(602,155){1}
\put(560,65){\circle*{10}}
\put(580,56){1}
\put(552,10){$I_2$}
\put(510,265){\circle*{10}}
\put(502,290){$I_6$}
\put(510,265){\line(1,-4){50}}
\put(506,155){3}

\put(1040,265){\circle*{10}}
\put(1032,290){$I_4$}
\put(1040,265){\line(-1,-4){50}}
\put(1030,155){2}
\put(990,65){\circle*{10}}
\put(1015,58){$\infty$}
\put(975,10){$I_2^{\ast}$}
\put(940,265){\circle*{10}}
\put(932,290){$I_4$}
\put(940,265){\line(1,-4){50}}
\put(936,155){2}

\put(1285,257){$\Phoch{1}$}
\put(1300,230){\vector(0,-1){125}}
\put(1285,57){$\Phoch{1}$}
\put(1325,160){$\pi_2$}

\end{picture}
\end{center}
\caption{A realization of \textbf{[2,2,4,4,6,6]}}
\label{Fig:2,2,4,4,6,6}
\end{figure}
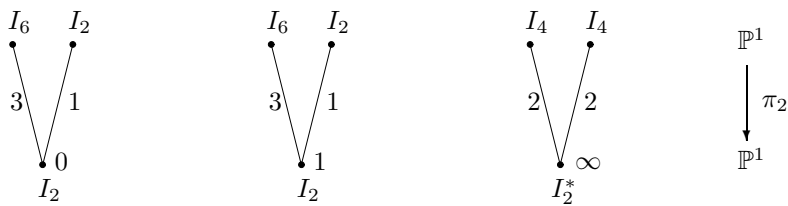
\vspace{0.2cm}

\section{The fibrations coming from degree 6 base changes}
\label{s:degree 6}

We consider the extremal rational elliptic surface
$\mathbf{X_{431}}$. By base change, it will give rise to 9
extremal elliptic K3 surfaces, 8 of which can be realized over
$\Q$.

We modify the Weierstrass equation of $\mathbf{X_{431}}$ given in
\cite{MP1} by exchanging 1 and $\infty$. It becomes
\[
y^2=x^3-3 \,(s-t)^3 (s-9 t)\,x-2 (s-t)^4 (s^2+18 s t-27 t^2).
\]
One finds the fibre of type $I_3$ above 0 and $I_1$ above $\infty$
while the $IV^{\ast}$ fibre sits above 1. We can resolve the
non-reduced fibre by the degree 6 base changes from Section
\ref{s:deflated}. Before the substitution, we can also permute 0
and $\infty$.

At first, we shall use $\mathbf{\pi_{3,4}}$. This leads to the
Weierstrass equation

\begin{small}
\begin{eqnarray*}
y^2 & = & x^3-3\,(-15 s^4 t^2+54 s^5 t+81 s^6+15 s^2 t^4-100 s^3 t^3+6 s t^5-t^6) (9 s^2+2 s t-t^2)\,x\\
&& \quad -2\,(19683 s^{12}+26244 s^{11} t+1458 s^{10} t^2+43740 s^9 t^3+25785 s^8 t^4-16776 s^7 t^5\\
&& \qquad\qquad -10108 s^6 t^6+3864 s^5 t^7+885 s^4 t^8-380 s^3
t^9-6 s^2 t^{10}+12 s t^{11}-t^{12}).
\end{eqnarray*}
\end{small}

This realizes \textbf{[1,2,3,3,3,12]}.

%\pi_3:
\begin{figure}[!ht]
\begin{center}
\setlength{\unitlength}{0.08mm}
\begin{picture}(1500,320)

\put(130,65){\circle*{10}}
\put(122,10){$I_3$}
\put(150,56){0}
\put(130,265){\circle*{10}}
\put(122,290){$I_3$}
\put(131,265){\line(0,-1){200}}
\put(105,155){1}
\put(50,265){\circle*{10}}
\put(35,290){$I_{12}$}
\put(50,265){\line(2,-5){80}}
\put(56,155){4}
\put(210,265){\circle*{10}}
\put(202,290){$I_3$}
\put(210,265){\line(-2,-5){80}}
\put(190,155){1}

\put(610,265){\circle*{10}} \put(602,290){$I_0$}
\put(610,265){\line(-1,-4){50}} \put(602,155){3}
\put(560,65){\circle*{10}} \put(580,56){1}
\put(545,10){$IV^{\ast}$} \put(510,265){\circle*{10}}
\put(502,290){$I_0$} \put(510,265){\line(1,-4){50}}
\put(506,155){3}
\put(990,65){\circle*{10}} \put(982,10){$I_1$}
\put(1010,58){$\infty$} \put(990,265){\circle*{10}}
\put(982,290){$I_2$} \put(991,265){\line(0,-1){200}}
\put(965,155){2} \put(910,265){\circle*{10}} \put(902,290){$I_3$}
\put(910,265){\line(2,-5){80}} \put(916,155){3}
\put(1070,265){\circle*{10}} \put(1062,290){$I_1$}
\put(1070,265){\line(-2,-5){80}} \put(1050,155){1}

\put(1285,257){$\Phoch{1}$}
\put(1300,230){\vector(0,-1){125}}
\put(1285,57){$\Phoch{1}$}
\put(1325,160){$\pi_{3,4}$}

\end{picture}
\end{center}
\caption{A realization of \textbf{[1,2,3,3,3,12]}}
\end{figure}
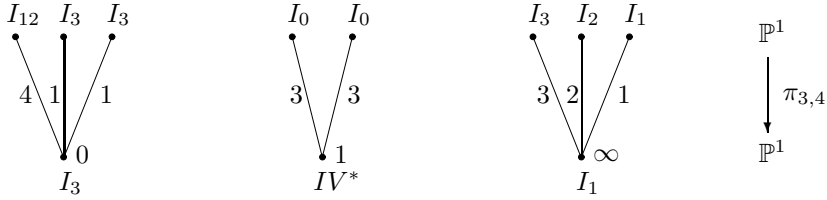
\vspace{0.2cm}

Permuting 0 and $\infty$ before the substitution gives a realization of \textbf{[1,1,3,4,6,9]}:
\begin{eqnarray*}
y^2 & = & x^3-3\,(-t^6+6 s t^5+15 s^2 t^4-100 s^3 t^3-1215 s^4 t^2+4374 s^5 t+6561 s^6)\\
&& \qquad\qquad (9 s^2+2 s t-t^2)\,x\\
&& \quad +2\, (14348907 s^{12}+19131876 s^{11} t+1062882 s^{10} t^2-4855140 s^9 t^3\\
&& \qquad\quad -185895 s^8 t^4+452952 s^7 t^5-7084 s^6 t^6-20328 s^5 t^7\\
&& \qquad\quad +3405 s^4 t^8-380 s^3 t^9-6 s^2 t^{10}+12 s t^{11}-t^{12}).
\end{eqnarray*}

%\pi_3:
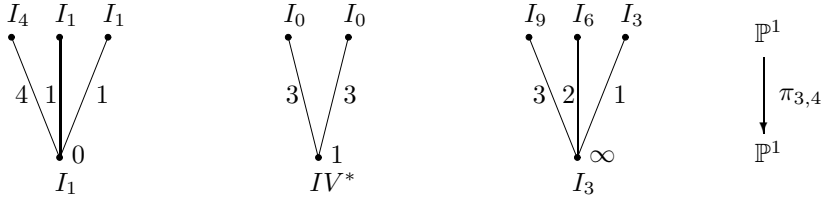
\begin{figure}[!ht]
\begin{center}
\setlength{\unitlength}{0.08mm}
\begin{picture}(1500,320)

\put(130,65){\circle*{10}}
\put(122,10){$I_1$}
\put(150,56){0}
\put(130,265){\circle*{10}}
\put(122,290){$I_1$}
\put(131,265){\line(0,-1){200}}
\put(105,155){1}
\put(50,265){\circle*{10}}
\put(42,290){$I_4$}
\put(50,265){\line(2,-5){80}}
\put(56,155){4}
\put(210,265){\circle*{10}}
\put(202,290){$I_1$}
\put(210,265){\line(-2,-5){80}}
\put(190,155){1}

\put(610,265){\circle*{10}}
\put(602,290){$I_0$}
\put(610,265){\line(-1,-4){50}}
\put(602,155){3}
\put(560,65){\circle*{10}}
\put(580,56){1}
\put(545,10){$IV^{\ast}$}
\put(510,265){\circle*{10}}
\put(502,290){$I_0$}
\put(510,265){\line(1,-4){50}}
\put(506,155){3}

\put(990,65){\circle*{10}}
\put(982,10){$I_3$}
\put(1010,58){$\infty$}
\put(990,265){\circle*{10}}
\put(982,290){$I_6$}
\put(991,265){\line(0,-1){200}}
\put(965,155){2}
\put(910,265){\circle*{10}}
\put(902,290){$I_9$}
\put(910,265){\line(2,-5){80}}
\put(916,155){3}
\put(1070,265){\circle*{10}}
\put(1062,290){$I_3$}
\put(1070,265){\line(-2,-5){80}}
\put(1050,155){1}

\put(1285,257){$\Phoch{1}$}
\put(1300,230){\vector(0,-1){125}}
\put(1285,57){$\Phoch{1}$}
\put(1325,160){$\pi_{3,4}$}

\end{picture}
\end{center}
\caption{A realization of \textbf{[1,1,3,4,6,9]}}
\end{figure}
\vspace{0.2cm}

We now turn to the substitutions by $\mathbf{\pi_{5,3}}$. These provide the following realizations of \textbf{[1,1,1,3,3,15]} as
\begin{eqnarray*}
y^2 & = & x^3-3\,(s^2-ts-t^2)(s^6-3s^5t+45t^3s^3-27t^5s-9t^6)\,x\\
&& \quad +2\,(-s^{12}+6s^{11}t-9s^{10}t^2+90s^9t^3-270t^4s^8-54t^5s^7+819t^6s^6\\
&& \qquad\quad
+54t^7s^5-810t^8s^4-270t^9s^3+243t^{10}s^2+162st^{11}+27t^{12})
\end{eqnarray*}

%Grad 6:
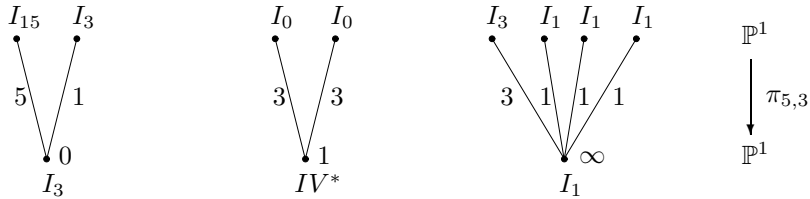
\begin{figure}[!ht]
\begin{center}
\setlength{\unitlength}{0.08mm}
\begin{picture}(1500,320)

\put(180,265){\circle*{10}}
\put(172,290){$I_3$}
\put(180,265){\line(-1,-4){50}}
\put(172,155){1}
\put(130,65){\circle*{10}}
\put(150,56){0}
\put(122,10){$I_3$}
\put(80,265){\circle*{10}}
\put(65,290){$I_{15}$}
\put(80,265){\line(1,-4){50}}
\put(76,155){5}

\put(610,265){\circle*{10}}
\put(602,290){$I_0$}
\put(610,265){\line(-1,-4){50}}
\put(602,155){3}
\put(560,65){\circle*{10}}
\put(580,56){1}
\put(542,10){$IV^{\ast}$}
\put(510,265){\circle*{10}}
\put(502,290){$I_0$}
\put(510,265){\line(1,-4){50}}
\put(506,155){3}

\put(1110,265){\circle*{10}}
\put(1102,290){$I_1$}
\put(1110,265){\line(-3,-5){120}}
\put(1070,155){1}
\put(1023,265){\circle*{10}}
\put(1015,290){$I_1$}
\put(1023,265){\line(-1,-6){33}}
\put(1012,155){1}
\put(990,65){\circle*{10}}
\put(1015,58){$\infty$}
\put(982,10){$I_1$}
\put(870,265){\circle*{10}}
\put(860,290){$I_3$}
\put(870,265){\line(3,-5){120}}
\put(884,155){3}
\put(957,265){\circle*{10}}
\put(949,290){$I_1$}
\put(957,265){\line(1,-6){33}}
\put(949,155){1}

\put(1285,257){$\Phoch{1}$}
\put(1300,230){\vector(0,-1){125}}
\put(1285,57){$\Phoch{1}$}
\put(1325,160){$\pi_{5,3}$}

\end{picture}
\end{center}
\caption{A realization of \textbf{[1,1,1,3,3,15]}}
\label{fig:15}
\end{figure}
\vspace{0.2cm}

and of \textbf{[1,3,3,3,5,9]} as
\begin{eqnarray*}
y^2 & = & x^3-3\-
,(s^2-s t-t^2)(5 s^3 t^3-3 s t^5-t^6+9 s^6-27 s^5 t)\,x\\
&& \quad +2\,(27 s^{12}-162 s^{11} t+243 s^{10} t^2+90 s^9 t^3-270 s^8 t^4-54 s^7 t^5\\
&& \qquad\quad +119 s^6 t^6+54 s^5 t^7+30 s^4 t^8+10 s^3 t^9-9 s^2 t^{10}-6 s t^{11}-t^{12}).
\end{eqnarray*}

%Grad 6:
\begin{figure}[!ht]
\begin{center}
\setlength{\unitlength}{0.08mm}
\begin{picture}(1500,320)

\put(180,265){\circle*{10}}
\put(172,290){$I_1$}
\put(180,265){\line(-1,-4){50}}
\put(172,155){1}
\put(130,65){\circle*{10}}
\put(150,56){0}
\put(122,10){$I_1$}
\put(80,265){\circle*{10}}
\put(72,290){$I_5$}
\put(80,265){\line(1,-4){50}}
\put(76,155){5}

\put(610,265){\circle*{10}}
\put(602,290){$I_0$}
\put(610,265){\line(-1,-4){50}}
\put(602,155){3}
\put(560,65){\circle*{10}}
\put(580,56){1}
\put(542,10){$IV^{\ast}$}
\put(510,265){\circle*{10}}
\put(502,290){$I_0$}
\put(510,265){\line(1,-4){50}}
\put(506,155){3}

\put(1110,265){\circle*{10}}
\put(1102,290){$I_3$}
\put(1110,265){\line(-3,-5){120}}
\put(1070,155){1}
\put(1023,265){\circle*{10}}
\put(1015,290){$I_3$}
\put(1023,265){\line(-1,-6){33}}
\put(1012,155){1}
\put(990,65){\circle*{10}}
\put(1015,58){$\infty$}
\put(982,10){$I_3$}
\put(870,265){\circle*{10}}
\put(860,290){$I_9$}
\put(870,265){\line(3,-5){120}}
\put(884,155){3}
\put(957,265){\circle*{10}}
\put(949,290){$I_3$}
\put(957,265){\line(1,-6){33}}
\put(949,155){1}

\put(1285,257){$\Phoch{1}$}
\put(1300,230){\vector(0,-1){125}}
\put(1285,57){$\Phoch{1}$}
\put(1325,160){$\pi_{5,3}$}

\end{picture}
\end{center}
\caption{A realization of \textbf{[1,3,3,3,5,9]}}
\end{figure}
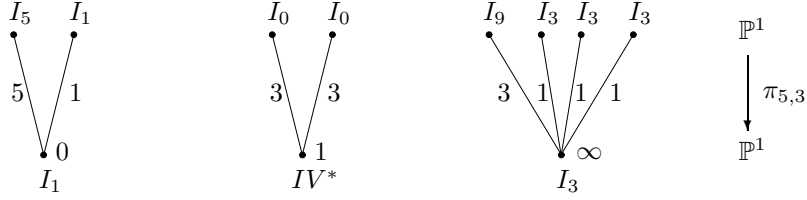
\vspace{0.2cm}

The substitution by $\mathbf{\pi_{5,2}}$ allows us to realize
\textbf{[1,1,2,2,3,15]} by virtue of the Weierstrass equation
\begin{eqnarray*}
y^2 & = & x^3-3\,(125 s^6-786 s^5 t+1575 s^4 t^2+1300 s^3 t^3+315 s^2 t^4+30 s t^5+t^6)\\
&& \qquad\quad (5 s^2+10 s t+t^2)\,x\\
&& \quad +2\,(15625 s^{12}+112986 s^{10} t^2-100500 s^{11} t-941300 s^9 t^3\\
&& \qquad\quad +1514175 s^8 t^4+3849240 s^7 t^5+2658380 s^6 t^6+912696 s^5 t^7\\
&& \qquad\quad +t^{12}+180375 s^4 t^8+21500 s^3 t^9+1530 s^2
t^{10}+60 s t^{11}).
\end{eqnarray*}

%Grad 6:
\begin{figure}[!ht]
\begin{center}
\setlength{\unitlength}{0.08mm}
\begin{picture}(1500,320)

\put(180,265){\circle*{10}}
\put(172,290){$I_3$}
\put(180,265){\line(-1,-4){50}}
\put(172,155){1}
\put(130,65){\circle*{10}}
\put(150,56){0}
\put(122,10){$I_3$}
\put(80,265){\circle*{10}}
\put(65,290){$I_{15}$}
\put(80,265){\line(1,-4){50}}
\put(76,155){5}

\put(610,265){\circle*{10}}
\put(602,290){$I_0$}
\put(610,265){\line(-1,-4){50}}
\put(602,155){3}
\put(560,65){\circle*{10}}
\put(580,56){1}
\put(542,10){$IV^{\ast}$}
\put(510,265){\circle*{10}}
\put(502,290){$I_0$}
\put(510,265){\line(1,-4){50}}
\put(506,155){3}

\put(1110,265){\circle*{10}}
\put(1102,290){$I_1$}
\put(1110,265){\line(-3,-5){120}}
\put(1070,155){1}
\put(1023,265){\circle*{10}}
\put(1015,290){$I_1$}
\put(1023,265){\line(-1,-6){33}}
\put(1012,155){1}
\put(990,65){\circle*{10}}
\put(1015,58){$\infty$}
\put(982,10){$I_1$}
\put(870,265){\circle*{10}}
\put(860,290){$I_2$}
\put(870,265){\line(3,-5){120}}
\put(884,155){2}
\put(957,265){\circle*{10}}
\put(949,290){$I_2$}
\put(957,265){\line(1,-6){33}}
\put(948,155){2}

\put(1285,257){$\Phoch{1}$}
\put(1300,230){\vector(0,-1){125}}
\put(1285,57){$\Phoch{1}$}
\put(1325,160){$\pi_{5,2}$}

\end{picture}
\end{center}
\caption{A realization of \textbf{[1,1,2,2,3,15]}}
\label{fig:20}
\end{figure}
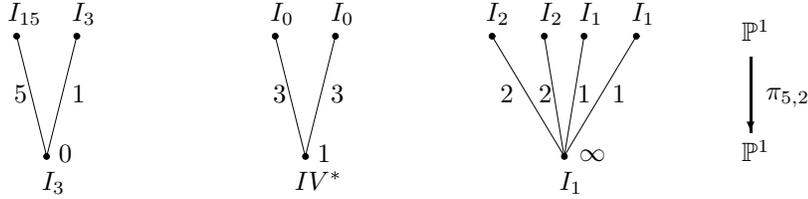
\vspace{0.2cm}

Then, the configuration \textbf{[1,3,3,5,6,6]} comes from
\begin{eqnarray*}
y^2 & = & x^3-3\,(125 s^6+14574 s^5 t+1575 s^4 t^2+1300 s^3 t^3+315 s^2 t^4+30 s t^5+t^6)\\
&& \qquad\quad (5 s^2+10 s t+t^2)\,x\\
&& \quad -2\,(15625 s^{12}-4132500 s^{11} t-48851622 s^{10} t^2-51744500 s^9 t^3\\
&& \qquad\quad -40418625 s^8 t^4-6311400 s^7 t^5+1690700 s^6 t^6+880440 s^5 t^7\\
&& \qquad\quad +180375 s^4 t^8+21500 s^3 t^9+1530 s^2 t^{10}+60 s t^{11}+t^{12}).
\end{eqnarray*}

%Grad 6:
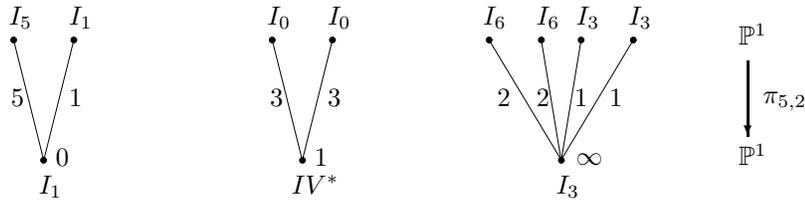
\begin{figure}[!ht]
\begin{center}
\setlength{\unitlength}{0.08mm}
\begin{picture}(1500,320)

\put(180,265){\circle*{10}}
\put(172,290){$I_1$}
\put(180,265){\line(-1,-4){50}}
\put(172,155){1}
\put(130,65){\circle*{10}}
\put(150,56){0}
\put(122,10){$I_1$}
\put(80,265){\circle*{10}}
\put(72,290){$I_5$}
\put(80,265){\line(1,-4){50}}
\put(76,155){5}

\put(610,265){\circle*{10}}
\put(602,290){$I_0$}
\put(610,265){\line(-1,-4){50}}
\put(602,155){3}
\put(560,65){\circle*{10}}
\put(580,56){1}
\put(542,10){$IV^{\ast}$}
\put(510,265){\circle*{10}}
\put(502,290){$I_0$}
\put(510,265){\line(1,-4){50}}
\put(506,155){3}

\put(1110,265){\circle*{10}}
\put(1102,290){$I_3$}
\put(1110,265){\line(-3,-5){120}}
\put(1070,155){1}
\put(1023,265){\circle*{10}}
\put(1015,290){$I_3$}
\put(1023,265){\line(-1,-6){33}}
\put(1012,155){1}
\put(990,65){\circle*{10}}
\put(1015,58){$\infty$}
\put(982,10){$I_3$}
\put(870,265){\circle*{10}}
\put(860,290){$I_6$}
\put(870,265){\line(3,-5){120}}
\put(884,155){2}
\put(957,265){\circle*{10}}
\put(949,290){$I_6$}
\put(957,265){\line(1,-6){33}}
\put(949,155){2}

\put(1285,257){$\Phoch{1}$}
\put(1300,230){\vector(0,-1){125}}
\put(1285,57){$\Phoch{1}$}
\put(1325,160){$\pi_{5,2}$}

\end{picture}
\end{center}
\caption{A realization of \textbf{[1,3,3,5,6,6]}}
\end{figure}
\vspace{0.2cm}

Finally, we can also substitute by $\mathbf{\pi_{4,3}}$ and
produce Weierstrass equations for \textbf{[1,1,1,3,6,12]} as
\begin{eqnarray*}
y^2 & = & x^3-3\,(-276s^4t^2+8s^6+96s^5t+416s^3t^3-186s^2t^4+24st^5-t^6)\\
&& \qquad\quad (2s^2+8st-t^2)\,x\\
&& \quad +2\,(11160s^8t^4+7392s^{10}t^2-15232s^9t^3-130176s^7t^5\\
&&\qquad\quad +220056s^6t^6-160416s^5t^7+54792s^4t^8+64s^{12}\\
&& \qquad\quad +1536s^{11}t-9760s^3t^9+948s^2t^{10}-48st^{11}+t^{12})
\end{eqnarray*}

%Grad 6:
\begin{figure}[!ht]
\begin{center}
\setlength{\unitlength}{0.08mm}
\begin{picture}(1500,320)

\put(180,265){\circle*{10}}
\put(172,290){$I_6$}
\put(180,265){\line(-1,-4){50}}
\put(172,155){2}
\put(130,65){\circle*{10}}
\put(150,56){0}
\put(122,10){$I_3$}
\put(80,265){\circle*{10}}
\put(65,290){$I_{12}$}
\put(80,265){\line(1,-4){50}}
\put(76,155){4}

\put(610,265){\circle*{10}}
\put(602,290){$I_0$}
\put(610,265){\line(-1,-4){50}}
\put(602,155){3}
\put(560,65){\circle*{10}}
\put(580,56){1}
\put(542,10){$IV^{\ast}$}
\put(510,265){\circle*{10}}
\put(502,290){$I_0$}
\put(510,265){\line(1,-4){50}}
\put(506,155){3}

\put(1110,265){\circle*{10}}
\put(1102,290){$I_1$}
\put(1110,265){\line(-3,-5){120}}
\put(1070,155){1}
\put(1023,265){\circle*{10}}
\put(1015,290){$I_1$}
\put(1023,265){\line(-1,-6){33}}
\put(1012,155){1}
\put(990,65){\circle*{10}}
\put(1015,58){$\infty$}
\put(982,10){$I_1$}
\put(870,265){\circle*{10}}
\put(860,290){$I_3$}
\put(870,265){\line(3,-5){120}}
\put(884,155){3}
\put(957,265){\circle*{10}}
\put(949,290){$I_1$}
\put(957,265){\line(1,-6){33}}
\put(949,155){1}

\put(1285,257){$\Phoch{1}$}
\put(1300,230){\vector(0,-1){125}}
\put(1285,57){$\Phoch{1}$}
\put(1325,160){$\pi_{4,3}$}

\end{picture}
\end{center}
\caption{A realization of \textbf{[1,1,1,3,6,12]}}
\end{figure}
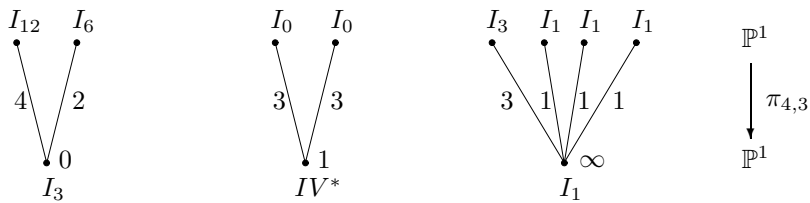
\vspace{0.2cm}

and for \textbf{[2,3,3,3,4,9]} as
\begin{eqnarray*}
y^2 & = & x^3-3\,(8s^6+96s^5t+6204s^4t^2+416s^3t^3-186s^2t^4+24st^5-t^6)\\
&& \qquad\quad (2s^2+8st-t^2)\,x\\
&& \quad -2\,(68400s^4t^8+64s^{12}-101472s^{10}t^2+1536s^{11}t\\
&& \qquad\quad +2751144s^6t^6-1321600s^9t^3-9460008s^8t^4-5791104s^7t^5\\
&& \qquad\quad -487008s^5t^7-9760s^3t^9+948s^2t^{10}-48st^{11}+t^{12}).
\end{eqnarray*}

%Grad 6:
\begin{figure}[!ht]
\begin{center}
\setlength{\unitlength}{0.08mm}
\begin{picture}(1500,320)

\put(180,265){\circle*{10}}
\put(172,290){$I_2$}
\put(180,265){\line(-1,-4){50}}
\put(172,155){2}
\put(130,65){\circle*{10}}
\put(150,56){0}
\put(122,10){$I_1$}
\put(80,265){\circle*{10}}
\put(72,290){$I_4$}
\put(80,265){\line(1,-4){50}}
\put(76,155){4}

\put(610,265){\circle*{10}}
\put(602,290){$I_0$}
\put(610,265){\line(-1,-4){50}}
\put(602,155){3}
\put(560,65){\circle*{10}}
\put(580,56){1}
\put(542,10){$IV^{\ast}$}
\put(510,265){\circle*{10}}
\put(502,290){$I_0$}
\put(510,265){\line(1,-4){50}}
\put(506,155){3}

\put(1110,265){\circle*{10}}
\put(1102,290){$I_3$}
\put(1110,265){\line(-3,-5){120}}
\put(1070,155){1}
\put(1023,265){\circle*{10}}
\put(1015,290){$I_3$}
\put(1023,265){\line(-1,-6){33}}
\put(1012,155){1}
\put(990,65){\circle*{10}}
\put(1015,58){$\infty$}
\put(982,10){$I_3$}
\put(870,265){\circle*{10}}
\put(860,290){$I_9$}
\put(870,265){\line(3,-5){120}}
\put(884,155){3}
\put(957,265){\circle*{10}}
\put(949,290){$I_3$}
\put(957,265){\line(1,-6){33}}
\put(949,155){1}

\put(1285,257){$\Phoch{1}$}
\put(1300,230){\vector(0,-1){125}}
\put(1285,57){$\Phoch{1}$}
\put(1325,160){$\pi_{4,3}$}

\end{picture}
\end{center}
\caption{A realization of \textbf{[2,3,3,3,4,9]}}
\end{figure}
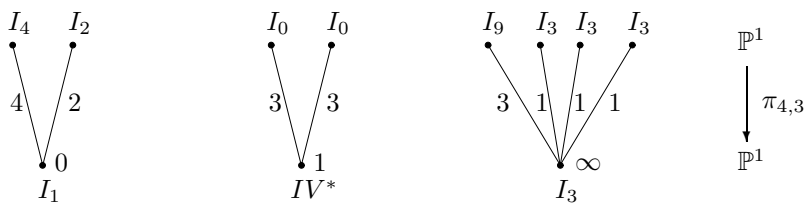
\vspace{0.2cm}

We remark that substitution by $\tilde{\pi}$ allows us to realize
the configuration \textbf{[1,2,3,3,6,9]} over the number field
$\Q(x^3+12x-12)$.

\section[The inflating base changes]{The inflating base changes\\ and the resulting fibrations}
\label{s:inflating}

We now turn to the other possibility to resolve non-reduced fibres
by a base change $\pi$. In this case we allow $\pi$ to be
inflating, i.e.~the pull-back $X$ may contain non-reduced fibres
of type $I_n^{\ast}$ for $n\geq 0$ apart from its semi-stable
fibres. The only additional assumption is that their number is
even. Via deflation, we can substitute these fibres by their
reduced relatives. The resulting semi-stable surface is a
quadratic twist of $X$.

One might hope to produce a number of new configurations by this
method. A close inspection, however, shows that none arise from
the extremal rational elliptic surfaces considered in the last two
sections. This approach is nevertheless useful, since it enables
us to work with the extremal rational elliptic surface
$\mathbf{X_{321}}$ as well. The reason is that the Hurwitz formula
is not violated if we choose the degree 8 base change $\pi$ of
$\Phoch{1}$ to have ramification index (2,2,2,2) at the cusp of
the $III^{\ast}$ fibre (instead of (4,4) or 8 before). The fibre
is thus replaced by four fibres of type $I_0^{\ast}$ which can
easily be twisted away.

We have to assume the other two cusps to have six pre-images in
total. Up to exchanging the cusps 0 and $\infty$, there are 13
such base changes which cannot be factored through an extremal
rational elliptic surface. In the following we concentrate on the
9 of these which can be defined over $\Q$. They result in 17
further extremal semi-stable elliptic K3 surfaces which arise by
pull-back from $\mathbf{X_{321}}$. The four base changes which are
not defined over $\Q$ will briefly be sketched at the end of this
section.

Remember that $\mathbf{X_{321}}$ has singular fibres of type
$III^{\ast}, I_2$ and $I_1$. We normalize the Weierstrass equation
given in \cite{MP1} such that the $III^{\ast}$ fibre sits above 1
and the $I_2$ and $I_1$ fibres above 0 and $\infty$, respectively:
\[
y^2=x^3-3\,(s-t)^3(s-4t)\,x-2\,(s-t)^5(s+8t).
\]
In the following, we investigate the degree 8 base changes of
$\Phoch{1}$ with ramification index (2,2,2,2) at 1 such that the
further cusps 0 and $\infty$ have six pre-images in total. We will
concentrate on those base changes which give rise to new
configurations of extremal elliptic K3 fibrations over $\Q$.
Without loss of generality, we assume that 0 does not have more
pre-images than $\infty$ (since we can exchange the cusps).

The three base changes which are totally ramified at 0 realize
configurations known from \cite{LY} or \cite{TY}. We now turn to
the base changes $\pi$ such that 0 has two pre-images. Throughout
we can exclude ramification of index (2,2,2,2) at both other cusps
since this would not allow ramification of high order at 0. We
list the base changes according to the ramification index at 0.

\vspace{0.2cm}

Let $\pi$ have ramification index \textbf{(7,1)} at 0. There are
four base changes to consider. The computations show that the base
change with ramification index (4,2,1,1) at $\infty$ can only be
defined over the Galois extensions $\Q(\sqrt{-7})$. It is given at
the end of this section. Here we give the remaining three base
changes over $\Q$:

The first base change has ramification index (5,1,1,1) at
$\infty$. It can be given as
\begin{eqnarray*}
\mathbf{\pi}:~~~ \Phoch{1}~ & \rightarrow & \quad\Phoch{1}\\
(s:t) & \mapsto & (s^7(s-4t):-4t^5(14s^3+14ts^2+20t^2s+25t^3))
\end{eqnarray*}
since
$s^7(s-4t)+4t^5(14s^3+14ts^2+20t^2s+25t^3)=(10t^4+4st^3+2s^2t^2-2s^3t-s^4)^2$.

Substituting $\pi$ into the normalized Weierstrass equation of
$\mathbf{X_{321}}$ gives, after deflation, a realization of
\textbf{[1,1,1,2,5,14]}:
\begin{eqnarray*}
y^2 & = & x^3-3\,(16s^8+32s^7t-112t^5s^3+56t^6s^2-40t^7s+25t^8)\,x\\
&& \quad -2\,(-5t^4+4t^3s-4s^2t^2+8s^3t+8s^4)\\
&& \qquad\quad (8s^8+16s^7t+112t^5s^3-56t^6s^2+40t^7s-25t^8).
\end{eqnarray*}

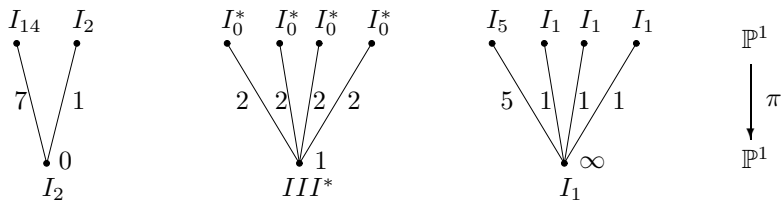
\begin{figure}[!ht]
\begin{center}
\setlength{\unitlength}{0.08mm}
\begin{picture}(1500,320)

\put(180,265){\circle*{10}}
\put(172,290){$I_2$}
\put(180,265){\line(-1,-4){50}}
\put(172,155){1}
\put(130,65){\circle*{10}}
\put(150,56){0}
\put(122,10){$I_2$}
\put(80,265){\circle*{10}}
\put(65,290){$I_{14}$}
\put(80,265){\line(1,-4){50}}
\put(76,155){7}

\put(670,265){\circle*{10}}
\put(662,290){$I_0^{\ast}$}
\put(670,265){\line(-3,-5){120}}
\put(630,155){2}
\put(583,265){\circle*{10}}
\put(575,290){$I_0^{\ast}$}
\put(583,265){\line(-1,-6){33}}
\put(572,155){2}
\put(550,65){\circle*{10}}
\put(575,56){1}
\put(521,10){$III^{\ast}$}
\put(430,265){\circle*{10}}
\put(420,290){$I_0^{\ast}$}
\put(430,265){\line(3,-5){120}}
\put(444,155){2}
\put(517,265){\circle*{10}}
\put(509,290){$I_0^{\ast}$}
\put(517,265){\line(1,-6){33}}
\put(509,155){2}

\put(1110,265){\circle*{10}}
\put(1102,290){$I_1$}
\put(1110,265){\line(-3,-5){120}}
\put(1070,155){1}
\put(1023,265){\circle*{10}}
\put(1015,290){$I_1$}
\put(1023,265){\line(-1,-6){33}}
\put(1012,155){1}
\put(990,65){\circle*{10}}
\put(1015,58){$\infty$}
\put(982,10){$I_1$}
\put(870,265){\circle*{10}}
\put(860,290){$I_5$}
\put(870,265){\line(3,-5){120}}
\put(884,155){5}
\put(957,265){\circle*{10}}
\put(949,290){$I_1$}
\put(957,265){\line(1,-6){33}}
\put(949,155){1}

\put(1285,257){$\Phoch{1}$}
\put(1300,230){\vector(0,-1){125}}
\put(1285,57){$\Phoch{1}$}
\put(1325,160){$\pi$}

\end{picture}
\end{center}
\caption{A realization of \textbf{[1,1,1,2,5,14]}}
\label{fig:35}
\end{figure}

Permuting the cusps 0 and $\infty$ before the substitution leads
to the following realization of \textbf{[1,2,2,2,7,10]}:
\begin{eqnarray*}
y^2 & = & x^3-3\,(-14t^5s^3+14t^6s^2-20t^7s+25t^8+s^8+4s^7t)\,x\\
&& \quad +(-10t^4+4t^3s-2s^2t^2+2s^3t+s^4)\\
&& \qquad\quad (14t^5s^3-14t^6s^2+20t^7s-25t^8+2s^8+8s^7t).
\end{eqnarray*}

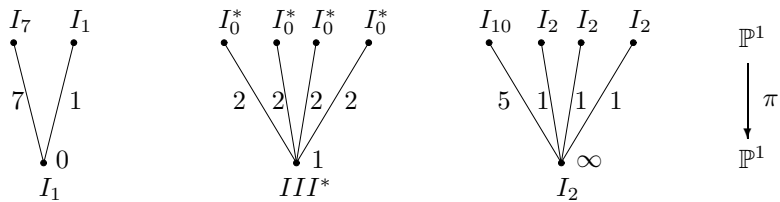
\begin{figure}[!ht]
\begin{center}
\setlength{\unitlength}{0.08mm}
\begin{picture}(1500,320)

\put(180,265){\circle*{10}}
\put(172,290){$I_1$}
\put(180,265){\line(-1,-4){50}}
\put(172,155){1}
\put(130,65){\circle*{10}}
\put(150,56){0}
\put(122,10){$I_1$}
\put(80,265){\circle*{10}}
\put(72,290){$I_7$}
\put(80,265){\line(1,-4){50}}
\put(76,155){7}

\put(670,265){\circle*{10}}
\put(662,290){$I_0^{\ast}$}
\put(670,265){\line(-3,-5){120}}
\put(630,155){2}
\put(583,265){\circle*{10}}
\put(575,290){$I_0^{\ast}$}
\put(583,265){\line(-1,-6){33}}
\put(572,155){2}
\put(550,65){\circle*{10}}
\put(575,56){1}
\put(521,10){$III^{\ast}$}
\put(430,265){\circle*{10}}
\put(420,290){$I_0^{\ast}$}
\put(430,265){\line(3,-5){120}}
\put(444,155){2}
\put(517,265){\circle*{10}}
\put(509,290){$I_0^{\ast}$}
\put(517,265){\line(1,-6){33}}
\put(509,155){2}

\put(1110,265){\circle*{10}}
\put(1102,290){$I_2$}
\put(1110,265){\line(-3,-5){120}}
\put(1070,155){1}
\put(1023,265){\circle*{10}}
\put(1015,290){$I_2$}
\put(1023,265){\line(-1,-6){33}}
\put(1012,155){1}
\put(990,65){\circle*{10}}
\put(1015,58){$\infty$}
\put(982,10){$I_2$}
\put(870,265){\circle*{10}}
\put(853,290){$I_{10}$}
\put(870,265){\line(3,-5){120}}
\put(884,155){5}
\put(957,265){\circle*{10}}
\put(949,290){$I_2$}
\put(957,265){\line(1,-6){33}}
\put(949,155){1}

\put(1285,257){$\Phoch{1}$}
\put(1300,230){\vector(0,-1){125}}
\put(1285,57){$\Phoch{1}$}
\put(1325,160){$\pi$}

\end{picture}
\end{center}
\caption{A realization of \textbf{[1,2,2,2,7,10]}}
\end{figure}
\vspace{0.2cm}

The second base change has ramification index (3,3,1,1) at
$\infty$. It can be chosen as
\begin{eqnarray*}
\mathbf{\pi}:\quad \Phoch{1}\; & \rightarrow & \quad\Phoch{1}\\
(s:t) & \mapsto & (1728s^7t:-(s^2-5st+t^2)^3(7s^2-13st+t^2))
\end{eqnarray*}
since
$1728s^7t+(s^2-5st+t^2)^3(7s^2-13st+t^2)=(t^4-14st^3+63s^2t^2-70s^3t-7s^4)^2$.

As pull-back of $\mathbf{X_{321}}$ via $\pi$ we realize the
constellations \textbf{[1,1,2,3,3,14]}:
\begin{scriptsize}
\begin{eqnarray*}
y^2 & = & x^3-3\,(49s^8-316s^7t+4018s^6t^2-8624s^5t^3+5915s^4t^4-1904s^3t^5+322s^2t^6-28st^7+t^8)\,x\\
&& \quad +2(49s^8-964s^7t+4018s^6t^2-8624s^5t^3+5915s^4t^4-1904s^3t^5+322s^2t^6-28st^7+t^8)\\
&& \qquad\quad (t^4-14st^3+63s^2t^2-70s^3t-7s^4)
\end{eqnarray*}
\end{scriptsize}

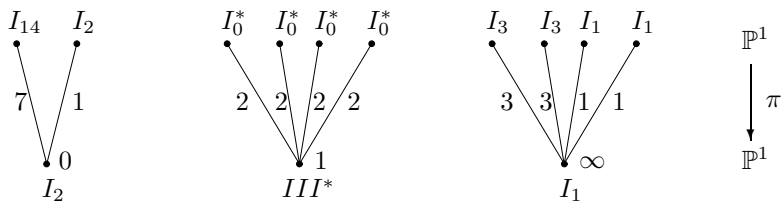
\begin{figure}[!ht]
\begin{center}
\setlength{\unitlength}{0.08mm}
\begin{picture}(1500,320)

\put(180,265){\circle*{10}}
\put(172,290){$I_2$}
\put(180,265){\line(-1,-4){50}}
\put(172,155){1}
\put(130,65){\circle*{10}}
\put(150,56){0}
\put(122,10){$I_2$}
\put(80,265){\circle*{10}}
\put(65,290){$I_{14}$}
\put(80,265){\line(1,-4){50}}
\put(76,155){7}

\put(670,265){\circle*{10}}
\put(662,290){$I_0^{\ast}$}
\put(670,265){\line(-3,-5){120}}
\put(630,155){2}
\put(583,265){\circle*{10}}
\put(575,290){$I_0^{\ast}$}
\put(583,265){\line(-1,-6){33}}
\put(572,155){2}
\put(550,65){\circle*{10}}
\put(575,56){1}
\put(521,10){$III^{\ast}$}
\put(430,265){\circle*{10}}
\put(420,290){$I_0^{\ast}$}
\put(430,265){\line(3,-5){120}}
\put(444,155){2}
\put(517,265){\circle*{10}}
\put(509,290){$I_0^{\ast}$}
\put(517,265){\line(1,-6){33}}
\put(509,155){2}

\put(1110,265){\circle*{10}}
\put(1102,290){$I_1$}
\put(1110,265){\line(-3,-5){120}}
\put(1070,155){1}
\put(1023,265){\circle*{10}}
\put(1015,290){$I_1$}
\put(1023,265){\line(-1,-6){33}}
\put(1012,155){1}
\put(990,65){\circle*{10}}
\put(1015,58){$\infty$}
\put(982,10){$I_1$}
\put(870,265){\circle*{10}}
\put(860,290){$I_3$}
\put(870,265){\line(3,-5){120}}
\put(884,155){3}
\put(957,265){\circle*{10}}
\put(949,290){$I_3$}
\put(957,265){\line(1,-6){33}}
\put(949,155){3}

\put(1285,257){$\Phoch{1}$}
\put(1300,230){\vector(0,-1){125}}
\put(1285,57){$\Phoch{1}$}
\put(1325,160){$\pi$}

\end{picture}
\end{center}
\caption{A realization of \textbf{[1,1,2,3,3,14]}}
\label{fig:84}
\end{figure}
\vspace{0.2cm}

and \textbf{[1,2,2,6,6,7]}:
\begin{scriptsize}
\begin{eqnarray*}
y^2 & = & x^3-3\,(49s^8+6164s^7t+4018s^6t^2-8624s^5t^3+5915s^4t^4-1904s^3t^5+322s^2t^6-28st^7+t^8)\,x\\
&& \quad -2\,(49s^8-14572s^7t+4018s^6t^2-8624s^5t^3+5915s^4t^4-1904s^3t^5+322s^2t^6-28st^7+t^8).\\
&& \qquad\quad (t^4-14st^3+63s^2t^2-70s^3t-7s^4).
\end{eqnarray*}
\end{scriptsize}

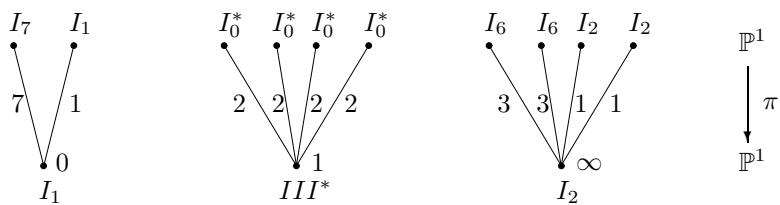
\begin{figure}[!ht]
\begin{center}
\setlength{\unitlength}{0.08mm}
\begin{picture}(1500,320)

\put(180,265){\circle*{10}}
\put(172,290){$I_1$}
\put(180,265){\line(-1,-4){50}}
\put(172,155){1}
\put(130,65){\circle*{10}}
\put(150,56){0}
\put(122,10){$I_1$}
\put(80,265){\circle*{10}}
\put(72,290){$I_7$}
\put(80,265){\line(1,-4){50}}
\put(76,155){7}

\put(670,265){\circle*{10}}
\put(662,290){$I_0^{\ast}$}
\put(670,265){\line(-3,-5){120}}
\put(630,155){2}
\put(583,265){\circle*{10}}
\put(575,290){$I_0^{\ast}$}
\put(583,265){\line(-1,-6){33}}
\put(572,155){2}
\put(550,65){\circle*{10}}
\put(575,56){1}
\put(521,10){$III^{\ast}$}
\put(430,265){\circle*{10}}
\put(420,290){$I_0^{\ast}$}
\put(430,265){\line(3,-5){120}}
\put(444,155){2}
\put(517,265){\circle*{10}}
\put(509,290){$I_0^{\ast}$}
\put(517,265){\line(1,-6){33}}
\put(509,155){2}

\put(1110,265){\circle*{10}}
\put(1102,290){$I_2$}
\put(1110,265){\line(-3,-5){120}}
\put(1070,155){1}
\put(1023,265){\circle*{10}}
\put(1015,290){$I_2$}
\put(1023,265){\line(-1,-6){33}}
\put(1012,155){1}
\put(990,65){\circle*{10}}
\put(1015,58){$\infty$}
\put(982,10){$I_2$}
\put(870,265){\circle*{10}}
\put(862,290){$I_6$}
\put(870,265){\line(3,-5){120}}
\put(884,155){3}
\put(957,265){\circle*{10}}
\put(949,290){$I_6$}
\put(957,265){\line(1,-6){33}}
\put(949,155){3}

\put(1285,257){$\Phoch{1}$}
\put(1300,230){\vector(0,-1){125}}
\put(1285,57){$\Phoch{1}$}
\put(1325,160){$\pi$}

\end{picture}
\end{center}
\caption{A realization of \textbf{[1,2,2,6,6,7]}}
\end{figure}
\vspace{0.2cm}

Recall that the pull-back surfaces inherit the group of sections
of the rational elliptic surfaces. As a consequence the above
fibration with configuration [1,1,2,3,3,14] necessarily has
Mordell-Weil group $\Z/(2)$. It differs from the surface with the
same configuration and $MW=(0)$, as obtained as a double sextic
over $\Q$ in \cite[p.55]{ATZ}. The underlying complex surfaces are not
isomorphic, since their discriminants differ.

The third base change has ramification index (3,2,2,1) at
$\infty$. We consider the map
\begin{eqnarray*}
\mathbf{\pi}:\quad \Phoch{1}\; & \rightarrow & \quad\Phoch{1}\\
(s:t) & \mapsto & (s^7(s+24t):16t^3(7s^2-14st+6t^2)^2(2s-t))
\end{eqnarray*}
with
$s^7(s+48t)-16t^3(7s^2-14st+6t^2)^2(2s-t)=(24t^4-80t^3s+72s^2t^2-12s^3t-s^4)^2$.

The pull-back via $\pi$ gives rise to the constellations \textbf{[1,2,2,2,3,14]}:
\begin{small}
\begin{eqnarray*}
y^2 & = & x^3-3\,(s^8+12s^7t-784t^3s^5+1764t^4s^4-1512t^5s^3+616t^6s^2-120t^7s+9t^8)\,x\\
&& \quad +(-s^8-12s^7t-1568t^3s^5+3528t^4s^4-3024t^5s^3+1232t^6s^2-240t^7s+18t^8)\\
&& \qquad\quad (-2s^4-12s^3t+36s^2t^2-20st^3+3t^4).
\end{eqnarray*}
\end{small}

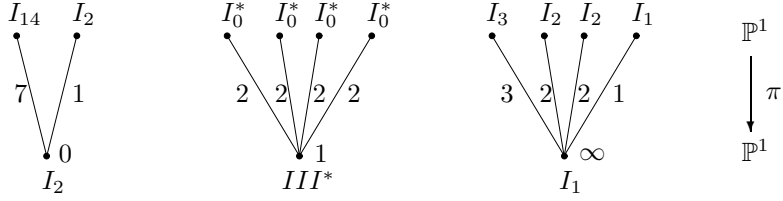
\begin{figure}[!ht]
\begin{center}
\setlength{\unitlength}{0.08mm}
\begin{picture}(1500,320)

\put(180,265){\circle*{10}}
\put(172,290){$I_2$}
\put(180,265){\line(-1,-4){50}}
\put(172,155){1}
\put(130,65){\circle*{10}}
\put(150,56){0}
\put(122,10){$I_2$}
\put(80,265){\circle*{10}}
\put(65,290){$I_{14}$}
\put(80,265){\line(1,-4){50}}
\put(76,155){7}

\put(670,265){\circle*{10}}
\put(662,290){$I_0^{\ast}$}
\put(670,265){\line(-3,-5){120}}
\put(630,155){2}
\put(583,265){\circle*{10}}
\put(575,290){$I_0^{\ast}$}
\put(583,265){\line(-1,-6){33}}
\put(572,155){2}
\put(550,65){\circle*{10}}
\put(575,56){1}
\put(521,10){$III^{\ast}$}
\put(430,265){\circle*{10}}
\put(420,290){$I_0^{\ast}$}
\put(430,265){\line(3,-5){120}}
\put(444,155){2}
\put(517,265){\circle*{10}}
\put(509,290){$I_0^{\ast}$}
\put(517,265){\line(1,-6){33}}
\put(509,155){2}

\put(1110,265){\circle*{10}}
\put(1102,290){$I_1$}
\put(1110,265){\line(-3,-5){120}}
\put(1070,155){1}
\put(1023,265){\circle*{10}}
\put(1015,290){$I_2$}
\put(1023,265){\line(-1,-6){33}}
\put(1012,155){2}
\put(990,65){\circle*{10}}
\put(1015,58){$\infty$}
\put(982,10){$I_1$}
\put(870,265){\circle*{10}}
\put(860,290){$I_3$}
\put(870,265){\line(3,-5){120}}
\put(884,155){3}
\put(957,265){\circle*{10}}
\put(949,290){$I_2$}
\put(957,265){\line(1,-6){33}}
\put(949,155){2}

\put(1285,257){$\Phoch{1}$}
\put(1300,230){\vector(0,-1){125}}
\put(1285,57){$\Phoch{1}$}
\put(1325,160){$\pi$}

\end{picture}
\end{center}
\caption{A realization of \textbf{[1,2,2,2,3,14]}}
\end{figure}
\vspace{0.2cm}

and \textbf{[1,2,4,4,6,7]}:
\begin{footnotesize}
\begin{eqnarray*}
y^2 & = & x^3-3\,(s^8-392s^5t^3+1764s^4t^4-3024s^3t^5+2464s^2t^6-960st^7+144t^8+24s^7 t)\,x\\
&& \quad +2\,(-24 t^4+80 s t^3-72 s^2 t^2+12 s^3 t+s^4)\\
&& \qquad\quad (196 s^5 t^3-882 s^4 t^4+1512 s^3 t^5-1232 s^2 t^6+480 s t^7-72 t^8+s^8+24 s^7 t).
\end{eqnarray*}
\end{footnotesize}

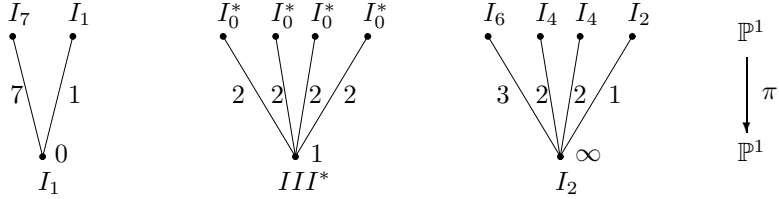
\begin{figure}[!ht]
\begin{center}
\setlength{\unitlength}{0.08mm}
\begin{picture}(1500,320)

\put(180,265){\circle*{10}}
\put(172,290){$I_1$}
\put(180,265){\line(-1,-4){50}}
\put(172,155){1}
\put(130,65){\circle*{10}}
\put(150,56){0}
\put(122,10){$I_1$}
\put(80,265){\circle*{10}}
\put(72,290){$I_7$}
\put(80,265){\line(1,-4){50}}
\put(76,155){7}

\put(670,265){\circle*{10}}
\put(662,290){$I_0^{\ast}$}
\put(670,265){\line(-3,-5){120}}
\put(630,155){2}
\put(583,265){\circle*{10}}
\put(575,290){$I_0^{\ast}$}
\put(583,265){\line(-1,-6){33}}
\put(572,155){2}
\put(550,65){\circle*{10}}
\put(575,56){1}
\put(521,10){$III^{\ast}$}
\put(430,265){\circle*{10}}
\put(420,290){$I_0^{\ast}$}
\put(430,265){\line(3,-5){120}}
\put(444,155){2}
\put(517,265){\circle*{10}}
\put(509,290){$I_0^{\ast}$}
\put(517,265){\line(1,-6){33}}
\put(509,155){2}

\put(1110,265){\circle*{10}}
\put(1102,290){$I_2$}
\put(1110,265){\line(-3,-5){120}}
\put(1070,155){1}
\put(1023,265){\circle*{10}}
\put(1015,290){$I_4$}
\put(1023,265){\line(-1,-6){33}}
\put(1012,155){2}
\put(990,65){\circle*{10}}
\put(1015,58){$\infty$}
\put(982,10){$I_2$}
\put(870,265){\circle*{10}}
\put(862,290){$I_6$}
\put(870,265){\line(3,-5){120}}
\put(884,155){3}
\put(957,265){\circle*{10}}
\put(949,290){$I_4$}
\put(957,265){\line(1,-6){33}}
\put(949,155){2}

\put(1285,257){$\Phoch{1}$}
\put(1300,230){\vector(0,-1){125}}
\put(1285,57){$\Phoch{1}$}
\put(1325,160){$\pi$}

\end{picture}
\end{center}
\caption{A realization of \textbf{[1,2,4,4,6,7]}}
\end{figure}
\vspace{0.2cm}

We turn to base changes with ramification index \textbf{(6,2)} at
0. Excluding all those ramification indices such that the
resulting configuration of singular fibres does not meet the
criteria of \cite{MP2}, the remaining base changes realize
configurations known from \cite{TY} or the previous sections. The
situation is exactly the same for ramification index
\textbf{(4,4)}. Hence, the only remaining ramification index at 0
is \textbf{(5,3)}. By the above considerations, we can exclude
ramification of index (3,3,1,1) at $\infty$. Then there are three
possibilities left:

For ramification index (5,1,1,1) at $\infty$, consider the map
\begin{eqnarray*}
\mathbf{\pi}:\quad \Phoch{1}\; & \rightarrow & \quad\Phoch{1}\\
(s:t) & \mapsto & (s^5(s-2t)^3:4t^5(6s^3-22s^2t-12st^2-9t^3))
\end{eqnarray*}
with
$s^5(s-2t)^3-4t^5(6s^3-22s^2t-12st^2-9t^3)=(6t^4+4st^3+6s^2t^2-6s^3t+s^4)^2$.

This base change realizes \textbf{[1,1,1,5,6,10]}:
\begin{small}
\begin{eqnarray*}
y^2 & = & x^3-3\,(16 s^8-96 s^7 t+192 s^6 t^2-128 s^5 t^3-48 s^3 t^5+88 s^2 t^6+24 s t^7+9 t^8)\,x\\
&& \quad -2\,(3 t^4+4 s t^3+12 s^2 t^2-24 s^3 t+8 s^4)\\
&& \qquad\quad (8 s^8-48 s^7 t+96 s^6 t^2-64 s^5 t^3+48 s^3 t^5-88 s^2 t^6-24 s t^7-9 t^8)
\end{eqnarray*}
\end{small}

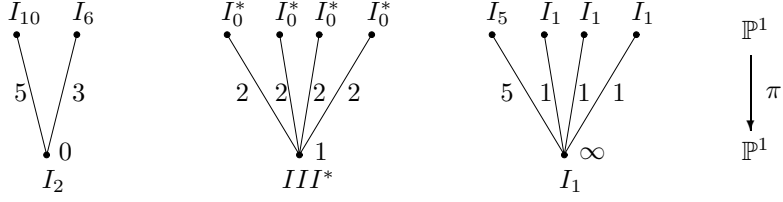
\begin{figure}[!ht]
\begin{center}
\setlength{\unitlength}{0.08mm}
\begin{picture}(1500,320)

\put(180,265){\circle*{10}}
\put(172,290){$I_6$}
\put(180,265){\line(-1,-4){50}}
\put(172,155){3}
\put(130,65){\circle*{10}}
\put(150,56){0}
\put(122,10){$I_2$}
\put(80,265){\circle*{10}}
\put(65,290){$I_{10}$}
\put(80,265){\line(1,-4){50}}
\put(76,155){5}

\put(670,265){\circle*{10}}
\put(662,290){$I_0^{\ast}$}
\put(670,265){\line(-3,-5){120}}
\put(630,155){2}
\put(583,265){\circle*{10}}
\put(575,290){$I_0^{\ast}$}
\put(583,265){\line(-1,-6){33}}
\put(572,155){2}
\put(550,65){\circle*{10}}
\put(575,56){1}
\put(521,10){$III^{\ast}$}
\put(430,265){\circle*{10}}
\put(420,290){$I_0^{\ast}$}
\put(430,265){\line(3,-5){120}}
\put(444,155){2}
\put(517,265){\circle*{10}}
\put(509,290){$I_0^{\ast}$}
\put(517,265){\line(1,-6){33}}
\put(509,155){2}

\put(1110,265){\circle*{10}}
\put(1102,290){$I_1$}
\put(1110,265){\line(-3,-5){120}}
\put(1070,155){1}
\put(1023,265){\circle*{10}}
\put(1015,290){$I_1$}
\put(1023,265){\line(-1,-6){33}}
\put(1012,155){1}
\put(990,65){\circle*{10}}
\put(1015,58){$\infty$}
\put(982,10){$I_1$}
\put(870,265){\circle*{10}}
\put(860,290){$I_5$}
\put(870,265){\line(3,-5){120}}
\put(884,155){5}
\put(957,265){\circle*{10}}
\put(949,290){$I_1$}
\put(957,265){\line(1,-6){33}}
\put(949,155){1}

\put(1285,257){$\Phoch{1}$}
\put(1300,230){\vector(0,-1){125}}
\put(1285,57){$\Phoch{1}$}
\put(1325,160){$\pi$}

\end{picture}
\end{center}
\caption{A realization of \textbf{[1,1,1,5,6,10]}}
\end{figure}
\vspace{0.2cm}

and after exchanging 0 and $\infty$ also \textbf{[2,2,2,3,5,10]}:
\begin{eqnarray*}
y^2 & = & x^3-3\,(s^8-6s^3t^5+22s^2t^6+12st^7+9t^8-12s^7t+48s^6t^2-64s^5t^3)\,x\\
&& \quad +(6t^4+4st^3+6s^2t^2-6s^3t+s^4)\\
&&
\quad\quad\;(6s^3t^5-22s^2t^6-12st^7-9t^8+2s^8-24s^7t+96s^6t^2-128s^5t^3).
\end{eqnarray*}

\begin{figure}[!ht]
\begin{center}
\setlength{\unitlength}{0.08mm}
\begin{picture}(1500,320)

\put(180,265){\circle*{10}}
\put(172,290){$I_3$}
\put(180,265){\line(-1,-4){50}}
\put(172,155){3}
\put(130,65){\circle*{10}}
\put(150,56){0}
\put(122,10){$I_1$}
\put(80,265){\circle*{10}}
\put(72,290){$I_5$}
\put(80,265){\line(1,-4){50}}
\put(76,155){5}

\put(670,265){\circle*{10}}
\put(662,290){$I_0^{\ast}$}
\put(670,265){\line(-3,-5){120}}
\put(630,155){2}
\put(583,265){\circle*{10}}
\put(575,290){$I_0^{\ast}$}
\put(583,265){\line(-1,-6){33}}
\put(572,155){2}
\put(550,65){\circle*{10}}
\put(575,56){1}
\put(521,10){$III^{\ast}$}
\put(430,265){\circle*{10}}
\put(420,290){$I_0^{\ast}$}
\put(430,265){\line(3,-5){120}}
\put(444,155){2}
\put(517,265){\circle*{10}}
\put(509,290){$I_0^{\ast}$}
\put(517,265){\line(1,-6){33}}
\put(509,155){2}

\put(1110,265){\circle*{10}}
\put(1102,290){$I_2$}
\put(1110,265){\line(-3,-5){120}}
\put(1070,155){1}
\put(1023,265){\circle*{10}}
\put(1015,290){$I_2$}
\put(1023,265){\line(-1,-6){33}}
\put(1012,155){1}
\put(990,65){\circle*{10}}
\put(1015,58){$\infty$}
\put(982,10){$I_2$}
\put(870,265){\circle*{10}}
\put(853,290){$I_{10}$}
\put(870,265){\line(3,-5){120}}
\put(884,155){5}
\put(957,265){\circle*{10}}
\put(949,290){$I_2$}
\put(957,265){\line(1,-6){33}}
\put(949,155){1}

\put(1285,257){$\Phoch{1}$}
\put(1300,230){\vector(0,-1){125}}
\put(1285,57){$\Phoch{1}$}
\put(1325,160){$\pi$}

\end{picture}
\end{center}
\caption{A realization of \textbf{[2,2,2,3,5,10]}}
\end{figure}
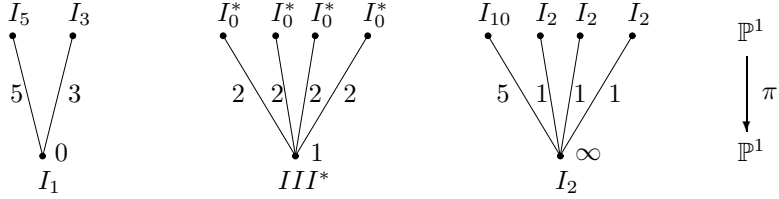
\vspace{0.2cm}

The ramification index (4,2,1,1) at $\infty$ can be obtained by the following base change:
\begin{eqnarray*}
\mathbf{\pi}:\quad \Phoch{1}\; & \rightarrow & \quad\Phoch{1}\\
(s:t) & \mapsto & (s^5(s-8t)^3:4t^4(3s+t)^2(5s^2-42st-9t^2))
\end{eqnarray*}
with
$s^5(s-8t)^3-4t^4(3s+t)^2(5s^2-42st-9t^2)=(s^4-12s^3t+24s^2t^2+32st^3+3t^4)^2$.

This realizes \textbf{[1,1,2,4,6,10]}:
\begin{eqnarray*}
y^2 & = & x^3-3\,(16 s^8-192 s^7 t+768 s^6 t^2-1024 s^5 t^3-720 s^4 t^4\\
&& \quad\qquad +2784 s^3 t^5+1312 s^2 t^6+192 s t^7+9 t^8)\,x\\
&&\quad -2\,(8 s^4-48 s^3 t+48 s^2 t^2+32 s t^3+3 t^4)\\
&& \quad\qquad (8 s^8-96 s^7 t+384 s^6 t^2-512 s^5 t^3+720 s^4 t^4\\
&& \qquad\qquad\; -2784 s^3 t^5-1312 s^2 t^6-192 s t^7-9 t^8).
\end{eqnarray*}

\begin{figure}[!ht]
\begin{center}
\setlength{\unitlength}{0.08mm}
\begin{picture}(1500,320)

\put(180,265){\circle*{10}}
\put(172,290){$I_6$}
\put(180,265){\line(-1,-4){50}}
\put(172,155){3}
\put(130,65){\circle*{10}}
\put(150,56){0}
\put(122,10){$I_2$}
\put(80,265){\circle*{10}}
\put(65,290){$I_{10}$}
\put(80,265){\line(1,-4){50}}
\put(76,155){5}

\put(670,265){\circle*{10}}
\put(662,290){$I_0^{\ast}$}
\put(670,265){\line(-3,-5){120}}
\put(630,155){2}
\put(583,265){\circle*{10}}
\put(575,290){$I_0^{\ast}$}
\put(583,265){\line(-1,-6){33}}
\put(572,155){2}
\put(550,65){\circle*{10}}
\put(575,56){1}
\put(521,10){$III^{\ast}$}
\put(430,265){\circle*{10}}
\put(420,290){$I_0^{\ast}$}
\put(430,265){\line(3,-5){120}}
\put(444,155){2}
\put(517,265){\circle*{10}}
\put(509,290){$I_0^{\ast}$}
\put(517,265){\line(1,-6){33}}
\put(509,155){2}

\put(1110,265){\circle*{10}}
\put(1102,290){$I_1$}
\put(1110,265){\line(-3,-5){120}}
\put(1070,155){1}
\put(1023,265){\circle*{10}}
\put(1015,290){$I_1$}
\put(1023,265){\line(-1,-6){33}}
\put(1012,155){1}
\put(990,65){\circle*{10}}
\put(1015,58){$\infty$}
\put(982,10){$I_1$}
\put(870,265){\circle*{10}}
\put(860,290){$I_4$}
\put(870,265){\line(3,-5){120}}
\put(884,155){4}
\put(957,265){\circle*{10}}
\put(949,290){$I_2$}
\put(957,265){\line(1,-6){33}}
\put(949,155){2}

\put(1285,257){$\Phoch{1}$}
\put(1300,230){\vector(0,-1){125}}
\put(1285,57){$\Phoch{1}$}
\put(1325,160){$\pi$}

\end{picture}
\end{center}
\caption{A realization of \textbf{[1,1,2,4,6,10]}}
\label{fig:5,3}
\end{figure}
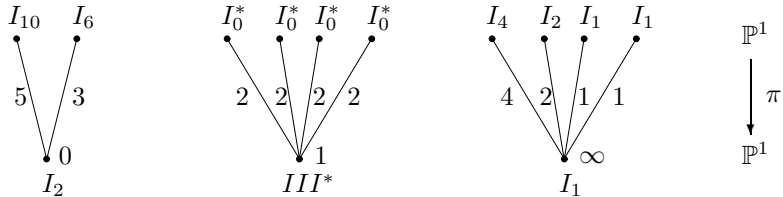
\vspace{0.2cm}

The permutation of 0 and $\infty$ leads to the constellation
\textbf{[2,2,3,4,5,8]}:
\begin{footnotesize}\begin{eqnarray*}
y^2 & = & x^3-3\,(s^8-24s^7t+192s^6t^2-512s^5t^3-45s^4t^4+348s^3t^5+328s^2t^6+96st^7+9t^8)\,x\\
&& \quad +(s^4-12s^3t+24s^2t^2+32st^3+6t^4)\\
&&\quad\quad\;
(2s^8-48s^7t+384s^6t^2-1024s^5t^3+45s^4t^4-348s^3t^5-328s^2t^6-96st^7-9t^8)
\end{eqnarray*}\end{footnotesize}

%Grad 8:
\begin{figure}[!ht]
\begin{center}
\setlength{\unitlength}{0.08mm}
\begin{picture}(1500,320)

\put(180,265){\circle*{10}}
\put(172,290){$I_3$}
\put(180,265){\line(-1,-4){50}}
\put(172,155){3}
\put(130,65){\circle*{10}}
\put(150,56){0}
\put(122,10){$I_1$}
\put(80,265){\circle*{10}}
\put(72,290){$I_5$}
\put(80,265){\line(1,-4){50}}
\put(76,155){5}

\put(670,265){\circle*{10}}
\put(662,290){$I_0^{\ast}$}
\put(670,265){\line(-3,-5){120}}
\put(630,155){2}
\put(583,265){\circle*{10}}
\put(575,290){$I_0^{\ast}$}
\put(583,265){\line(-1,-6){33}}
\put(572,155){2}
\put(550,65){\circle*{10}}
\put(575,56){1}
\put(521,10){$III^{\ast}$}
\put(430,265){\circle*{10}}
\put(420,290){$I_0^{\ast}$}
\put(430,265){\line(3,-5){120}}
\put(444,155){2}
\put(517,265){\circle*{10}}
\put(509,290){$I_0^{\ast}$}
\put(517,265){\line(1,-6){33}}
\put(509,155){2}

\put(1110,265){\circle*{10}}
\put(1102,290){$I_2$}
\put(1110,265){\line(-3,-5){120}}
\put(1070,155){1}
\put(1023,265){\circle*{10}}
\put(1015,290){$I_2$}
\put(1023,265){\line(-1,-6){33}}
\put(1012,155){1}
\put(990,65){\circle*{10}}
\put(1015,58){$\infty$}
\put(982,10){$I_2$}
\put(870,265){\circle*{10}}
\put(862,290){$I_8$}
\put(870,265){\line(3,-5){120}}
\put(884,155){4}
\put(957,265){\circle*{10}}
\put(949,290){$I_4$}
\put(957,265){\line(1,-6){33}}
\put(949,155){2}

\put(1285,257){$\Phoch{1}$}
\put(1300,230){\vector(0,-1){125}}
\put(1285,57){$\Phoch{1}$}
\put(1325,160){$\pi$}

\end{picture}
\end{center}
\caption{A realization of \textbf{[2,2,3,4,5,8]}}
\label{fig:5,3'}
\end{figure}
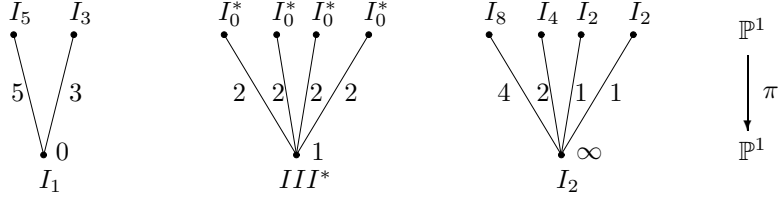
\vspace{0.2cm}

The final possible ramification index at $\infty$ is (3,2,2,1).
This is encoded in the following base change:
\begin{eqnarray*}
\mathbf{\pi}:\quad \Phoch{1}\; & \rightarrow & \quad\Phoch{1}\\
(s:t) & \mapsto & (4s^5(9s-4t)^3:-t^3(4s+t)(10s^2-6st+t^2)^2)
\end{eqnarray*}
with
$4s^5(9s-4t)^3+t^3(4s+t)(10s^2-6st+t^2)^2=(54s^4-36s^3t+4s^2t^2+4st^3-t^4)^2$.

This base change enables us to realize \textbf{[1,2,2,3,6,10]}:
\begin{footnotesize}
\begin{eqnarray*}
y^2 & = & x^3-3\,(9s^8-36s^7t+48s^6t^2+112s^5t^3-380s^4t^4+312s^3t^5+72s^2t^6-216st^7+81t^8)\,x\\
&& \quad +(-6s^4+12s^3t-4s^2t^2-12st^3+9t^4)\\
&& \quad\quad\; (-9s^8+36s^7t-48s^6t^2+288s^5t^3-760s^4t^4+624s^3t^5+144s^2t^6-432st^7+162t^8)
\end{eqnarray*}
\end{footnotesize}

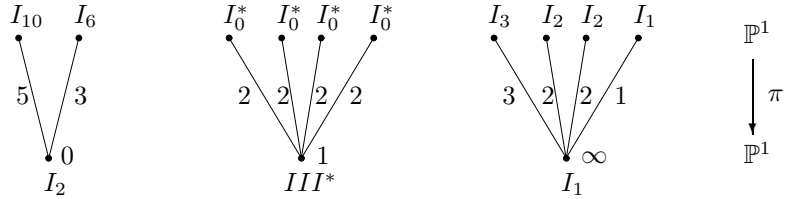
\begin{figure}[!ht]
\begin{center}
\setlength{\unitlength}{0.08mm}
\begin{picture}(1500,320)

\put(180,265){\circle*{10}}
\put(172,290){$I_6$}
\put(180,265){\line(-1,-4){50}}
\put(172,155){3}
\put(130,65){\circle*{10}}
\put(150,56){0}
\put(122,10){$I_2$}
\put(80,265){\circle*{10}}
\put(65,290){$I_{10}$}
\put(80,265){\line(1,-4){50}}
\put(76,155){5}

\put(670,265){\circle*{10}}
\put(662,290){$I_0^{\ast}$}
\put(670,265){\line(-3,-5){120}}
\put(630,155){2}
\put(583,265){\circle*{10}}
\put(575,290){$I_0^{\ast}$}
\put(583,265){\line(-1,-6){33}}
\put(572,155){2}
\put(550,65){\circle*{10}}
\put(575,56){1}
\put(521,10){$III^{\ast}$}
\put(430,265){\circle*{10}}
\put(420,290){$I_0^{\ast}$}
\put(430,265){\line(3,-5){120}}
\put(444,155){2}
\put(517,265){\circle*{10}}
\put(509,290){$I_0^{\ast}$}
\put(517,265){\line(1,-6){33}}
\put(509,155){2}

\put(1110,265){\circle*{10}}
\put(1102,290){$I_1$}
\put(1110,265){\line(-3,-5){120}}
\put(1070,155){1}
\put(1023,265){\circle*{10}}
\put(1015,290){$I_2$}
\put(1023,265){\line(-1,-6){33}}
\put(1012,155){2}
\put(990,65){\circle*{10}}
\put(1015,58){$\infty$}
\put(982,10){$I_1$}
\put(870,265){\circle*{10}}
\put(860,290){$I_3$}
\put(870,265){\line(3,-5){120}}
\put(884,155){3}
\put(957,265){\circle*{10}}
\put(949,290){$I_2$}
\put(957,265){\line(1,-6){33}}
\put(949,155){2}

\put(1285,257){$\Phoch{1}$}
\put(1300,230){\vector(0,-1){125}}
\put(1285,57){$\Phoch{1}$}
\put(1325,160){$\pi$}

\end{picture}
\end{center}
\caption{A realization of \textbf{[1,2,2,3,6,10]}}
\end{figure}
\vspace{0.2cm}

and furthermore \textbf{[2,3,4,4,5,6]}:
\begin{eqnarray*}
y^2 & = & x^3-3\,(-208s^5t^3-380s^4t^4+312s^3t^5+72s^2t^6-216st^7\\
&& \qquad\qquad+81t^8+144s^8-576s^7t+768s^6t^2)\,x\\
&& \quad -2\,(-6s^4+12s^3t-4s^2t^2-12st^3+9t^4)\\
&& \quad\qquad (816s^5t^3-380s^4t^4+312s^3t^5+72s^2t^6-216st^7\\
&& \qquad\qquad+81t^8-288s^8+1152s^7t-1536s^6t^2).
\end{eqnarray*}

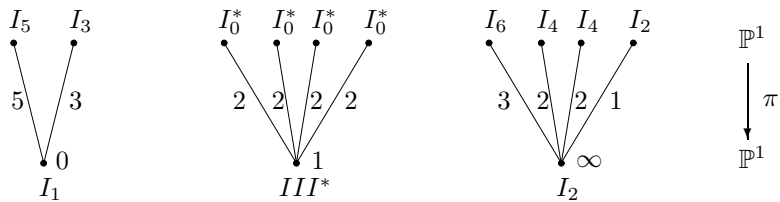
\begin{figure}[!ht]
\begin{center}
\setlength{\unitlength}{0.08mm}
\begin{picture}(1500,320)

\put(180,265){\circle*{10}}
\put(172,290){$I_3$}
\put(180,265){\line(-1,-4){50}}
\put(172,155){3}
\put(130,65){\circle*{10}}
\put(150,56){0}
\put(122,10){$I_1$}
\put(80,265){\circle*{10}}
\put(72,290){$I_5$}
\put(80,265){\line(1,-4){50}}
\put(76,155){5}

\put(670,265){\circle*{10}}
\put(662,290){$I_0^{\ast}$}
\put(670,265){\line(-3,-5){120}}
\put(630,155){2}
\put(583,265){\circle*{10}}
\put(575,290){$I_0^{\ast}$}
\put(583,265){\line(-1,-6){33}}
\put(572,155){2}
\put(550,65){\circle*{10}}
\put(575,56){1}
\put(521,10){$III^{\ast}$}
\put(430,265){\circle*{10}}
\put(420,290){$I_0^{\ast}$}
\put(430,265){\line(3,-5){120}}
\put(444,155){2}
\put(517,265){\circle*{10}}
\put(509,290){$I_0^{\ast}$}
\put(517,265){\line(1,-6){33}}
\put(509,155){2}

\put(1110,265){\circle*{10}}
\put(1102,290){$I_2$}
\put(1110,265){\line(-3,-5){120}}
\put(1070,155){1}
\put(1023,265){\circle*{10}}
\put(1015,290){$I_4$}
\put(1023,265){\line(-1,-6){33}}
\put(1012,155){2}
\put(990,65){\circle*{10}}
\put(1015,58){$\infty$}
\put(982,10){$I_2$}
\put(870,265){\circle*{10}}
\put(862,290){$I_6$}
\put(870,265){\line(3,-5){120}}
\put(884,155){3}
\put(957,265){\circle*{10}}
\put(949,290){$I_4$}
\put(957,265){\line(1,-6){33}}
\put(949,155){2}

\put(1285,257){$\Phoch{1}$}
\put(1300,230){\vector(0,-1){125}}
\put(1285,57){$\Phoch{1}$}
\put(1325,160){$\pi$}

\end{picture}
\end{center}
\caption{A realization of \textbf{[2,3,4,4,5,6]}}
\end{figure}
\vspace{0.2cm}

We shall now consider those base changes such that both cusps, 0
and $\infty$, have three pre-images. As a starting point, consider
ramification of index \textbf{(6,1,1)} at 0. There are three
respective ramification indices at $\infty$ such that $\pi$ does
not admit a factorization. However, (4,2,2) cannot occur, since
the resulting fibration does not exist by \cite{MP2}. The base
change with index (4,3,1) can only be defined over the number
field $\Q(\sqrt{-3})$. It is given at the end of this section. The
remaining base change with ramification index (5,2,1) at $\infty$
can be defined by
\begin{eqnarray*}
\mathbf{\pi}:\quad \Phoch{1}\; & \rightarrow & \quad\Phoch{1}\\
(s:t) & \mapsto & (4s^6(9s^2+24st+70 t^2):t^5(14s-5t)^2(4s-t))
\end{eqnarray*}
with
$4s^6(9s^2+24st+70t^2)-t^5(14s-5t)^2(4s-t)=(5t^4-24t^3s+18s^2t^2+8s^3t+6s^4)^2.$

The pull-back surface by $\pi$ has the singular fibres \textbf{[1,2,2,2,5,12]}:
\begin{eqnarray*}
y^2 & = & x^3-3\,(9s^8+24s^7t+70s^6t^2-784t^5s^3+756t^6s^2-240t^7s+25t^8)\,x\\
&& \quad +(50t^82-9s^8-24s^7t-70s^6t^2-1568t^5s^3+1512t^6s^2-480t^7s)\\
&& \quad\quad\; (5t^4-24st^3+18s^2t^2+8s^3t+6s^4)
\end{eqnarray*}

%degree 8:
\begin{figure}[!ht]
\begin{center}
\setlength{\unitlength}{0.08mm}
\begin{picture}(1500,320)

\put(130,65){\circle*{10}}
\put(122,10){$I_2$}
\put(150,56){0}
\put(130,265){\circle*{10}}
\put(122,290){$I_2$}
\put(131,265){\line(0,-1){200}}
\put(105,155){1}
\put(50,265){\circle*{10}}
\put(35,290){$I_{12}$}
\put(50,265){\line(2,-5){80}}
\put(56,155){6}
\put(210,265){\circle*{10}}
\put(202,290){$I_2$}
\put(210,265){\line(-2,-5){80}}
\put(190,155){1}

\put(670,265){\circle*{10}}
\put(662,290){$I_0^{\ast}$}
\put(670,265){\line(-3,-5){120}}
\put(630,155){2}
\put(583,265){\circle*{10}}
\put(575,290){$I_0^{\ast}$}
\put(583,265){\line(-1,-6){33}}
\put(572,155){2}
\put(550,65){\circle*{10}}
\put(575,56){1}
\put(521,10){$III^{\ast}$}
\put(430,265){\circle*{10}}
\put(420,290){$I_0^{\ast}$}
\put(430,265){\line(3,-5){120}}
\put(444,155){2}
\put(517,265){\circle*{10}}
\put(509,290){$I_0^{\ast}$}
\put(517,265){\line(1,-6){33}}
\put(509,155){2}

\put(990,65){\circle*{10}}
\put(982,10){$I_1$}
\put(1000,58){$\infty$}
\put(990,265){\circle*{10}}
\put(982,290){$I_2$}
\put(991,265){\line(0,-1){200}}
\put(965,155){2}
\put(910,265){\circle*{10}}
\put(902,290){$I_5$}
\put(910,265){\line(2,-5){80}}
\put(916,155){5}
\put(1070,265){\circle*{10}}
\put(1062,290){$I_1$}
\put(1070,265){\line(-2,-5){80}}
\put(1050,155){1}

\put(1285,257){$\Phoch{1}$}
\put(1300,230){\vector(0,-1){125}}
\put(1285,57){$\Phoch{1}$}
\put(1325,160){$\pi$}

\end{picture}
\end{center}
\caption{A realization of \textbf{[1,2,2,2,5,12]}}
\label{fig:5,2,1}
\end{figure}
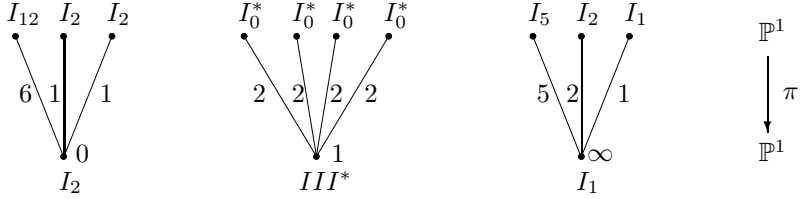
\vspace{0.2cm}

Permuting 0 and $\infty$ gives rise to the constellation
\textbf{[1,1,2,4,6,10]}. This was already realized in the
preceding paragraphs where 0 was assumed to have only two
pre-images (Fig. \ref{fig:5,3}). Although the fibrations differ,
the underlying complex surfaces are isomorphic. This can be
derived from lattice theory using the discriminant form
(cf.~\cite{SZ} for details).

We now come to base changes with ramification index
\textbf{(5,2,1)} at 0. Again, the other cusp $\infty$ cannot have
ramification index (4,2,2). Furthermore, the base change with
ramification index (5,2,1) at $\infty$ can only be defined over
the number field $\Q(7x^3+19x^2+16x+8)$. It is given at the end of
this section. Here, we only construct the two remaining base
changes with ramification index (4,3,1) or (3,3,2) at $\infty$.
For the first, consider the map
\begin{eqnarray*}
\mathbf{\pi}:\quad \Phoch{1}\; & \rightarrow & \quad\Phoch{1}\\
(s:t) & \mapsto & (2^8s^5(24s-7t)^2(4s+3t):-t^4(7s-t)^3(15s-t))
\end{eqnarray*}
with
$2^8s^5(24s-7t)^2(4s+3t)+t^4(7s-t)^3(15s-t)=(t^4-18t^3s+69s^2t^2-32s^3t-384s^4)^2.$

This enables us to realize the configurations \textbf{[1,2,3,4,4,10]}:
\begin{eqnarray*}
y^2 & = & x^3-3\,(36864 s^8+6144 s^7 t-12992 s^6 t^2+2352 t^3 s^5\\
&& \qquad\qquad +5145 t^4 s^4-2548 t^5 s^3+462 t^6 s^2-36 t^7 s+t^8)\,x\\
&& \quad -2\,(-t^4+18 s t^3-69 s^2 t^2+32 s^3 t+384 s^4)\\
&& \qquad\quad (18432 s^8+3072 s^7 t-6496 s^6 t^2+1176 t^3 s^5\\
&& \qquad\qquad -5145 t^4 s^4+2548 t^5 s^3-462 t^6 s^2+36 t^7 s-t^8)
\end{eqnarray*}

%\pi_3:
\begin{figure}[!ht]
\begin{center}
\setlength{\unitlength}{0.08mm}
\begin{picture}(1500,320)

\put(130,65){\circle*{10}}
\put(122,10){$I_2$}
\put(150,56){0}
\put(130,265){\circle*{10}}
\put(122,290){$I_4$}
\put(131,265){\line(0,-1){200}}
\put(105,155){2}
\put(50,265){\circle*{10}}
\put(35,290){$I_{10}$}
\put(50,265){\line(2,-5){80}}
\put(56,155){5}
\put(210,265){\circle*{10}}
\put(202,290){$I_2$}
\put(210,265){\line(-2,-5){80}}
\put(190,155){1}

\put(670,265){\circle*{10}}
\put(662,290){$I_0^{\ast}$}
\put(670,265){\line(-3,-5){120}}
\put(630,155){2}
\put(583,265){\circle*{10}}
\put(575,290){$I_0^{\ast}$}
\put(583,265){\line(-1,-6){33}}
\put(572,155){2}
\put(550,65){\circle*{10}}
\put(575,56){1}
\put(521,10){$III^{\ast}$}
\put(430,265){\circle*{10}}
\put(420,290){$I_0^{\ast}$}
\put(430,265){\line(3,-5){120}}
\put(444,155){2}
\put(517,265){\circle*{10}}
\put(509,290){$I_0^{\ast}$}
\put(517,265){\line(1,-6){33}}
\put(509,155){2}

\put(990,65){\circle*{10}}
\put(982,10){$I_1$}
\put(1000,58){$\infty$}
\put(990,265){\circle*{10}}
\put(982,290){$I_3$}
\put(991,265){\line(0,-1){200}}
\put(965,155){3}
\put(910,265){\circle*{10}}
\put(902,290){$I_4$}
\put(910,265){\line(2,-5){80}}
\put(916,155){4}
\put(1070,265){\circle*{10}}
\put(1062,290){$I_1$}
\put(1070,265){\line(-2,-5){80}}
\put(1050,155){1}

\put(1285,257){$\Phoch{1}$}
\put(1300,230){\vector(0,-1){125}}
\put(1285,57){$\Phoch{1}$}
\put(1325,160){$\pi$}

\end{picture}
\end{center}
\caption{A realization of \textbf{[1,2,3,4,4,10]}}
\end{figure}
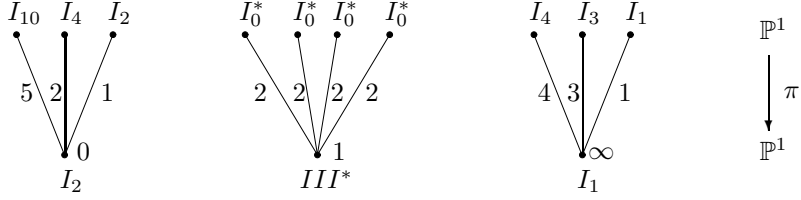
\vspace{0.2cm}

and \textbf{[1,2,2,5,6,8]}:
\begin{eqnarray*}
y^2 & = & x^3-3\,(5145 t^4 s^4-2548 t^5 s^3+462 t^6 s^2-36 t^7 s+t^8\\
&& \qquad\qquad +589824 s^8+98304 s^7 t-207872 s^6 t^2+37632 t^3 s^5)\,x\\
&& \quad +2\,(-t^4+18 s t^3-69 s^2 t^2+32 s^3 t+384 s^4)\\
&& \qquad\quad (-5145 t^4 s^4+2548 t^5 s^3-462 t^6 s^2+36 t^7 s-t^8\\
&& \qquad\qquad +1179648 s^8+196608 s^7 t-415744 s^6 t^2+75264 t^3 s^5).
\end{eqnarray*}

%\pi_3:
\begin{figure}[!ht]
\begin{center}
\setlength{\unitlength}{0.08mm}
\begin{picture}(1500,320)

\put(130,65){\circle*{10}}
\put(122,10){$I_1$}
\put(150,56){0}
\put(130,265){\circle*{10}}
\put(122,290){$I_2$}
\put(131,265){\line(0,-1){200}}
\put(105,155){2}
\put(50,265){\circle*{10}}
\put(42,290){$I_5$}
\put(50,265){\line(2,-5){80}}
\put(56,155){5}
\put(210,265){\circle*{10}}
\put(202,290){$I_1$}
\put(210,265){\line(-2,-5){80}}
\put(190,155){1}

\put(670,265){\circle*{10}}
\put(662,290){$I_0^{\ast}$}
\put(670,265){\line(-3,-5){120}}
\put(630,155){2}
\put(583,265){\circle*{10}}
\put(575,290){$I_0^{\ast}$}
\put(583,265){\line(-1,-6){33}}
\put(572,155){2}
\put(550,65){\circle*{10}}
\put(575,56){1}
\put(521,10){$III^{\ast}$}
\put(430,265){\circle*{10}}
\put(420,290){$I_0^{\ast}$}
\put(430,265){\line(3,-5){120}}
\put(444,155){2}
\put(517,265){\circle*{10}}
\put(509,290){$I_0^{\ast}$}
\put(517,265){\line(1,-6){33}}
\put(509,155){2}

\put(990,65){\circle*{10}}
\put(982,10){$I_2$}
\put(1000,58){$\infty$}
\put(990,265){\circle*{10}}
\put(982,290){$I_6$}
\put(991,265){\line(0,-1){200}}
\put(965,155){3}
\put(910,265){\circle*{10}}
\put(902,290){$I_8$}
\put(910,265){\line(2,-5){80}}
\put(916,155){4}
\put(1070,265){\circle*{10}}
\put(1062,290){$I_2$}
\put(1070,265){\line(-2,-5){80}}
\put(1050,155){1}

\put(1285,257){$\Phoch{1}$}
\put(1300,230){\vector(0,-1){125}}
\put(1285,57){$\Phoch{1}$}
\put(1325,160){$\pi$}

\end{picture}
\end{center}
\caption{A realization of \textbf{[1,2,2,5,6,8]}}
\end{figure}
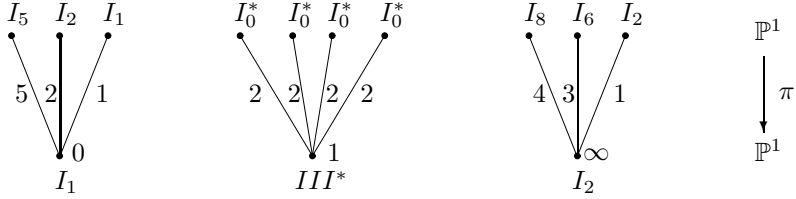
\vspace{0.2cm}

The second base change with ramification index (3,3,2) at $\infty$
can be constructed in the following way:
\begin{eqnarray*}
\mathbf{\pi}:\quad \Phoch{1}\; & \rightarrow & \quad\Phoch{1}\\
(s:t) & \mapsto & (9s^5(s+6t)^2(9s+4t):-4t^2(10s^2+24st+9t^2)^3)
\end{eqnarray*}
with
$9s^5(s+6t)^2(9s+4t)+4t^2(10s^2+24st+9t^2)^3=(9s^4+56s^3t+234s^2t^2+216st^3+54t^4)^2.$

It enables us to realize the configurations
\textbf{[2,2,3,3,4,10]}:
\begin{eqnarray*}
y^2 & = & x^3-3\,(1296s^8+8064s^7t+77392s^6t^2+232992s^5t^3+319680s^4t^4\\
&& \qquad\qquad +214272s^3t^5+71928s^2t^6+11664st^7+729t^8)\,x\\
&& \quad -2\,(27t^4+216st^3+468s^2t^2+224s^3t+72s^4)\\
&& \qquad\quad (648s^8+4032s^7t-57304s^6t^2-229104s^5t^3-319680s^4t^4\\
&& \qquad\qquad -214272s^3t^5-71928s^2t^6-11664st^7-729t^8)
\end{eqnarray*}

%\pi_3:
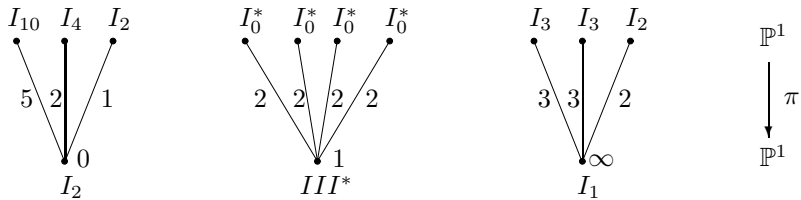
\begin{figure}[!ht]
\begin{center}
\setlength{\unitlength}{0.08mm}
\begin{picture}(1500,320)

\put(130,65){\circle*{10}}
\put(122,10){$I_2$}
\put(150,56){0}
\put(130,265){\circle*{10}}
\put(122,290){$I_4$}
\put(131,265){\line(0,-1){200}}
\put(105,155){2}
\put(50,265){\circle*{10}}
\put(35,290){$I_{10}$}
\put(50,265){\line(2,-5){80}}
\put(56,155){5}
\put(210,265){\circle*{10}}
\put(202,290){$I_2$}
\put(210,265){\line(-2,-5){80}}
\put(190,155){1}

\put(670,265){\circle*{10}}
\put(662,290){$I_0^{\ast}$}
\put(670,265){\line(-3,-5){120}}
\put(630,155){2}
\put(583,265){\circle*{10}}
\put(575,290){$I_0^{\ast}$}
\put(583,265){\line(-1,-6){33}}
\put(572,155){2}
\put(550,65){\circle*{10}}
\put(575,56){1}
\put(521,10){$III^{\ast}$}
\put(430,265){\circle*{10}}
\put(420,290){$I_0^{\ast}$}
\put(430,265){\line(3,-5){120}}
\put(444,155){2}
\put(517,265){\circle*{10}}
\put(509,290){$I_0^{\ast}$}
\put(517,265){\line(1,-6){33}}
\put(509,155){2}

\put(990,65){\circle*{10}}
\put(982,10){$I_1$}
\put(1000,58){$\infty$}
\put(990,265){\circle*{10}}
\put(982,290){$I_3$}
\put(991,265){\line(0,-1){200}}
\put(965,155){3}
\put(910,265){\circle*{10}}
\put(902,290){$I_3$}
\put(910,265){\line(2,-5){80}}
\put(916,155){3}
\put(1070,265){\circle*{10}}
\put(1062,290){$I_2$}
\put(1070,265){\line(-2,-5){80}}
\put(1050,155){2}

\put(1285,257){$\Phoch{1}$}
\put(1300,230){\vector(0,-1){125}}
\put(1285,57){$\Phoch{1}$}
\put(1325,160){$\pi$}

\end{picture}
\end{center}
\caption{A realization of \textbf{[2,2,3,3,4,10]}}
\label{fig:40}
\end{figure}
\vspace{0.2cm}

and \textbf{[1,2,4,5,6,6]}:
\begin{eqnarray*}
y^2 & = & x^3-3\,(4348s^6t^2+8496s^5t^3+19980t^4s^4+26784t^5s^3\\
&& \qquad\qquad +17982t^6s^2+5832t^7s+729t^8+81s^8+1008s^7t)\,x\\
&& \quad +\,(54t^4+216st^3+234s^2t^2+56s^3t+9s^4)\\
&& \qquad\quad (5696s^6t^2-4608s^5t^3-19980t^4s^4-26784t^5s^3\\
&& \qquad\qquad -17982t^6s^2-5832t^7s-729t^8+162s^8+2016s^7t).
\end{eqnarray*}

%\pi_3:
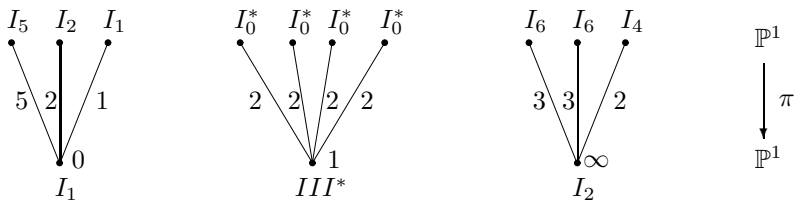
\begin{figure}[!ht]
\begin{center}
\setlength{\unitlength}{0.08mm}
\begin{picture}(1500,320)

\put(130,65){\circle*{10}}
\put(122,10){$I_1$}
\put(150,56){0}
\put(130,265){\circle*{10}}
\put(122,290){$I_2$}
\put(131,265){\line(0,-1){200}}
\put(105,155){2}
\put(50,265){\circle*{10}}
\put(42,290){$I_5$}
\put(50,265){\line(2,-5){80}}
\put(56,155){5}
\put(210,265){\circle*{10}}
\put(202,290){$I_1$}
\put(210,265){\line(-2,-5){80}}
\put(190,155){1}

\put(670,265){\circle*{10}}
\put(662,290){$I_0^{\ast}$}
\put(670,265){\line(-3,-5){120}}
\put(630,155){2}
\put(583,265){\circle*{10}}
\put(575,290){$I_0^{\ast}$}
\put(583,265){\line(-1,-6){33}}
\put(572,155){2}
\put(550,65){\circle*{10}}
\put(575,56){1}
\put(521,10){$III^{\ast}$}
\put(430,265){\circle*{10}}
\put(420,290){$I_0^{\ast}$}
\put(430,265){\line(3,-5){120}}
\put(444,155){2}
\put(517,265){\circle*{10}}
\put(509,290){$I_0^{\ast}$}
\put(517,265){\line(1,-6){33}}
\put(509,155){2}

\put(990,65){\circle*{10}}
\put(982,10){$I_2$}
\put(1000,58){$\infty$}
\put(990,265){\circle*{10}}
\put(982,290){$I_6$}
\put(991,265){\line(0,-1){200}}
\put(965,155){3}
\put(910,265){\circle*{10}}
\put(902,290){$I_6$}
\put(910,265){\line(2,-5){80}}
\put(916,155){3}
\put(1070,265){\circle*{10}}
\put(1062,290){$I_4$}
\put(1070,265){\line(-2,-5){80}}
\put(1050,155){2}

\put(1285,257){$\Phoch{1}$}
\put(1300,230){\vector(0,-1){125}}
\put(1285,57){$\Phoch{1}$}
\put(1325,160){$\pi$}

\end{picture}
\end{center}
\caption{A realization of \textbf{[1,2,4,5,6,6]}}
\end{figure}
\vspace{0.2cm}

Eventually, we come to the remaining base changes with
ramification index \textbf{(4,3,1)}, \textbf{(4,2,2)} or
\textbf{(3,3,2)} at 0. All but one of them either cannot exist or
give rise to known configurations of singular fibres. The
remaining base change has ramification index (4,3,1) at both
cusps, 0 and $\infty$. It can only be defined over the quadratic
extension $\Q(\sqrt{2})$ of $\Q$. We give it below.

\vspace{0.1cm}

We conclude this section by collecting the four base changes which
appeared above, but are only defined over an extension of $\Q$.
They are presented in the order of appearance in the course of
this section:

\begin{itemize}

\item A base change with ramification indices (7,1) and (4,2,1,1)
at 0 and $\infty$: Let $v$ be a solution of $2x^2-7x+28$. Then the
base change can be given by
$\pi((s:t))=(s^7(s+2vt):(10633/4v-2401)/40\,t^4(s-t)^2(s^2+(6-2v)st/5+(3v-14)t^2/10))$.
It realizes the configurations \textbf{[1,1,2,2,4,14]} and
\textbf{[1,2,2,4,7,8]} over $\Q(\sqrt{-7})$. This is the minimal
field of definition not only for the fibrations, but for the
surfaces. This follows from \cite[Prop. 13.1]{S}: The surfaces
have discriminant $d=56$ and $224$, so
$\Q(\sqrt{-d})=\Q(\sqrt{-14})$. This field has exponent 4.

\item The base change with ramification indices (6,1,1) and
(4,3,1) can be defined over $\Q(\sqrt{-3})$. If $v$ is a solution
to $3x^2-3x+7$, then we have $\pi((s:t))=
(s^6(s^2+4vst-(19v+14)/5\,t^2):(1763v-259)/20\,t^4(s-t)^3(s-(4v-7)/5\,t))$.
On the one hand we can thereby produce the configuration
\textbf{[1,1,2,6,6,8]}. The field of definition of this fibration
is the same as for the double sextic in \cite[p.~313]{P}. The
other pull-back has the configuration \textbf{[1,2,2,3,4,12]}. A
fibration with this configuration has already been obtained over
$\Q$ in Figure \ref{Fig:1,2,2,3,4,12}.

We will now show that these two fibrations with configuration
\textbf{[1,2,2,3,4,12]} have different Mordell-Weil groups: For
each fibration, consider the pull-back of a primitive section of
the basic rational elliptic surface $\mathbf{X_{321}}$
resp.~$\mathbf{X_{141}}$. This section can be directly computed in
terms of the components of the singular fibres which it meets.
Comparing this shape to \cite[Remark 0.4 (5)]{ATZ}, we conclude
that in both cases the induced section generates the Mordell-Weil
group of the pull-back. In particular, $MW$ of basic surface and
pull-back coincide. Since $MW(\mathbf{X_{321}})=\Z/2$ and
$MW(\mathbf{X_{141}})=\Z/4$, we obtain the claim.

The above fibrations over $\Q(\sqrt{-3})$ provide the first
examples where a fibration cannot be defined over $\Q$, although the underlying surface can. To see
this, we use the intersection form on the transcendental lattice
and the classification of \cite{SZ}.

By \cite{SZ}, the configuration \textbf{[1,1,2,6,6,8]} implies the
intersection form $diag(6,24)$. The complex K3 surface with this
form admits another elliptic fibration, No.~144 in the notation of
\cite{SZ}. Table \ref{T:extr-5} derives this fibration over $\Q$.
Similarly, the configuration \textbf{[1,2,2,3,4,12]} with
$MW={\Z}/2$ as above implies intersection form $12 I$. This
coincides with the intersection form of the \textbf{[2,2,4,4,6,6]}
configuration. This fibration was realized over $\Q$ in Figure
\ref{Fig:2,2,4,4,6,6}.

\item A base change with ramification index (5,2,1) at both cusps,
0 and $\infty$: Let $v$ be a zero of $7x^3+19x^2+16x+8$. Then
consider
$\pi((s:t))=(167s^5(s-2t)^2(s+4(v+1)\,t):-(15v^2+55v+52)\,t^5
(4s+(3v^2+3v-4)\,t)^2(8s-(7v^2+15v+4)\,t))$. The pull-back gives
rise to the configuration \textbf{[1,2,2,4,5,10]} over the
extension $\Q(x^3-75x+5150)$.

\item The final base change of this section has ramification index
(4,3,1) at 0 and $\infty$. It is defined over $\Q(\sqrt{2})$. Let
$v$ be a root of $7x^2+8x+2$. Then consider
$\pi((s:t))=(2^47^3s^4(s-t)^3(s+(8v+3)\,t):(8v+3)\,t^4(14s+(9v+4)\,t)^3(2s-(5v+4)\,t))$.
This map leads to the extremal K3 surface with singular fibres
\textbf{[1,2,3,4,6,8]}.

\end{itemize}

\vspace{0.1cm}

Together with \cite{TY}, the previous three sections exhaust the
extremal semi-stable elliptic K3 surfaces which can be realized as
non-general pull-back of rational elliptic surfaces. Here the term
"general" refers to the general pull-back construction involving
the induced $J$-map and the rational elliptic surface with
singular fibres $III^{\ast},II,I_1$. As a base change of degree 24
with very restricted ramification, this essentially makes no use
of the basic rational elliptic surface. A general solution for
this case has recently been announced by Beukers and Montanus
\cite{BM}. By construction, their fields of definition necessarily
coincide with ours in the overlapping cases.

\section{The non-semi-stable fibrations}
\label{s:non-semi-stable}

This section is devoted to a brief analysis of the non-semi-stable
extremal elliptic K3 fibrations. We will determine all fibrations
which can be derived from rational elliptic surfaces. The
treatment is significantly simplified by the fact that every such
surface has a non-reduced fibre. This follows from
\cite[Thm.~1.2]{Kl}. Moreover, every K3 fibration with more than
one non-reduced fibre is easily transformed into a rational
elliptic surface by the deflation process described in Section
\ref{s:basics}. Such an extremal K3 surface necessarily has three
or four cusps. In the case of three cusps, we find the K3 surfaces
directly in \cite{Ng}. If the K3 has four cusps, the corresponding
rational elliptic surfaces is given in \cite{H}. All but one of
the corresponding rational surfaces can be uniquely defined over
$\Q$ up to $\C$-isomorphism. Except for three cases which are
specified below, the deflation is also defined over $\Q$.

Table \ref{Table:rat*} collects the extremal K3 fibrations with
three or four cusps. The numbering refers to the classification in
\cite{SZ} and will be employed throughout this section. Each
fibration can be derived from a rational elliptic surface by
manipulating the Weierstrass equation. In other words, we reobtain
the rational elliptic surface by deflation.

\begin{table}[ht!]
\begin{small}
\begin{tabular}{|c|c|}
\hline No. & Config. \\
\hline
\hline 113 & 5,5,1*,1* \\
\hline 121 & 2,8,1*,1* \\
\hline 124 & 1,9,1*,1* \\
\hline 136 & 2*,2*,2* \\
\hline 137 & 4,4,2*,2* \\
\hline 153 & 3,6,1*,2* \\
\hline 154 & 1,8,1*,2* \\
\hline 155 & 3,3,3*,3* \\
\hline 167 & 2,6,1*,3* \\
\hline 168 & 1,6,2*,3* \\
\hline 169 & 2,2,4*,4* \\
\hline 177 & 1*,1*,4* \\
\hline 178 & 2,4,2*,4* \\
\hline 179 & 1,1,5*,5* \\
\hline 187 & 1,5,1*,5* \\
\hline 195 & 2,3,1*,6* \\
\hline 196 & 1,3,2*,6* \\
\hline 197 & 1,2,3*,6* \\
\hline 205 & 1,2,1*,8* \\
\hline
\end{tabular}
\begin{tabular}{|c|c|}
\hline No. & Config. \\
\hline
\hline 206 & 1,1,2*,8* \\
\hline 209 & 1,1,1*,9* \\
\hline 219 & IV*,IV*,IV*\\
\hline 220 & 4,4,IV*,IV* \\
\hline 222 & 2,4,IV*,IV*\\
\hline 226 & 1,7,IV*,IV* \\
\hline 243 & 4,5,1*,IV* \\
\hline 244 & 2,7,1*,IV* \\
\hline 245 & 1,8,1*,IV* \\
\hline 246 & 2*,IV*,IV* \\
\hline 247 & 3,5,2*,IV* \\
\hline 248 & 1,7,2*,IV* \\
\hline 249 & 2,5,3*.IV* \\
\hline 250 & 1*,3*,IV* \\
\hline 251 & 1,5,4*,IV* \\
\hline 252 & 2,3,5*,IV* \\
\hline 253 & 1,4,5*,IV* \\
\hline 254 & 1,2,7,IV* \\
\hline
\end{tabular}
\begin{tabular}{|c|c|}
\hline No. & Config. \\
\hline
\hline 255 & 1,1,8*,IV* \\
\hline 256 & 3,3,III*,III* \\
\hline 257 & 2,4,III*,III*\\
\hline 258 & 1,5,III*,III* \\
\hline 279 & 0*,III*,III* \\
\hline 280 & 3,5,1*,III* \\
\hline 281 & 2,6,1*,III* \\
\hline 282 & 1,7,1*,III* \\
\hline 283 & 3,4,2*,III* \\
\hline 284 & 1,6,2*,III* \\
\hline 285 & 1*,2*,III* \\
\hline 286 & 2,4,3*,III* \\
\hline 287 & 1,5,3*,III* \\
\hline 288 & 2,3,4*,III* \\
\hline 289 & 1,3,5*,III* \\
\hline 290 & 1,2,6*,III* \\
\hline 291 & 1,1,7*,III* \\
\hline 292 & 3,4,III*,IV* \\
\hline
%\phantom{292 & 3,4,III*,IV* \\ \hline}
\end{tabular}
\begin{tabular}{|c|c|}
\hline No. & Config. \\
\hline
\hline 293 & 2,5,III*,IV* \\
\hline 294 & 1,6,III*,IV* \\
\hline 295 & 1*,III*,IV* \\
\hline 296 & 2,2 II*,II*\\
\hline 297 & IV,II*,II*\\
\hline 313 & 1*,1*,II* \\
\hline 314 & 2,5,1*,II* \\
\hline 315 & 1,6,1*,II* \\
\hline 316 & 3,3,2*,II* \\
\hline 317 & 1,5,2*,II* \\
\hline 318 & 2,3,3*,II* \\
\hline 319 & 1,2,5*,II* \\
\hline 320 & 1,1,6*,II* \\
\hline 321 & 2,4,II*,IV* \\
\hline 322 & 1,5,II*,IV* \\
\hline 323 & 0*,II*,IV* \\
\hline 324 & 2,3,II*,III* \\
\hline 325 & 1,4,II*,III* \\
\hline
\end{tabular}
\vspace{0.2cm}
\caption{The extremal K3 fibrations with three or four cusps}
\label{Table:rat*}
\end{small}
\end{table}

There are four K3 fibrations in the above table which cannot be
defined over $\Q$ by this approach. For three of them, the basic
rational elliptic surfaces nevertheless is defined over $\Q$.
However, since it has cusps which are conjugate in some quadratic
field, certain manipulations are not defined over $\Q$. As a
result, for the fibrations with No.~187, No.~245 and No.~282,
fields of definition can only be given as $\Q(\sqrt{5}),
\Q(\sqrt{-2})$, and $\Q(\sqrt{-7})$, respectively. The rational
elliptic surface which corresponds to No.~294 can itself only be
defined over $\Q(\sqrt{-3})$.

For these four fibrations, we shall briefly discuss whether the
underlying K3 surface can be defined over $\Q$. This is possible
for No.~187, 245, and 294: Again we use the intersection form on
the transcendental lattice, as determined in \cite{SZ}. This gives
respective isomorphisms with the following surfaces: No.~259 from
Table \ref{T:extr-5}, No.~33 from Figure \ref{Fig:2,2,2,4,6,8},
and No.~84 from Figure \ref{Fig:1,2,2,3,4,12} (which has
$MW={\Z}/4$). On the other hand, the surface admitting
configuration No.~282 cannot be defined over $\Q$ by
\cite[Prop.~13.1]{S}. Hence the field of definition
$\Q(\sqrt{-7})$ is minimal.

%Trivial: Add stars to a rational elliptic surface with four or less cusps (all surfaces defined over $\Q$ unless otherwise noted):
%(cf.~\cite{Ng} for the examples with only three cusps)
%No.~113, 124, 136, 137,153, 154, 155, 167, 168, 169, 177, 178, 179, 187 (only over $\Q(\sqrt{5})$), 195, 196, 197, 205, 206, 209, 219, %220, 222, 226, 243-258 (No.~245 only over $\Q(\sqrt{-2})$), 279-297 (No.~282 only over $\Q(\sqrt{-7})$, No.~294 only over %$\Q(\sqrt{-3})$), 313-325.

%\vspace{0.1cm}

We now come to the extremal K3 fibrations which are still missing
with respect to the classification in \cite{SZ}. These have
exactly one non-reduced fibre. We will derive half of them from
rational elliptic surfaces. We will apply base change and further
make use of the "transfer of *" as explained in \cite{M2}.
Essentially this just moves the * from one fibre to another (a
priori not necessarily singular). This can be achieved by
translating the common factor of the polynomials $A$ and $B$ in
the Weierstrass equation. We will take base changes of degree 2 to
6 into account, depending on the singular fibres of the basic
rational elliptic surface. For each degree, we are going to
exploit one example in more detail. The remaining fibrations will
only be sketched very roughly.

\vspace{0.1cm}

\textbf{Degree 2:} These base changes involve the rational
elliptic surfaces with four singular fibres, which have one fibre
of type $III$. For example, take the surface $Y$ with singular
fibres $I_1,I_3,I_5,III$. According to \cite{H}, this surface and
all its cusps can be defined over $\Q$. Consider a quadratic base
change $\pi$ of $\Phoch{1}$ which is ramified at the cusp of the
$III$-fibre and at one further cusp. The pull-back of $Y$ via
$\pi$ is a K3 surface over $\Q$ with five semi-stable singular
fibres and one of type $I_0^{\ast}$. Transferring the * gives rise
to three extremal K3 surfaces with one non-reduced and four
semi-stable fibres. For instance, the configuration [2,3,3,5,5,0*]
can be transformed to \textbf{[3,3,5,5,2*]} (No.~138),
\textbf{[2,3,5,5,3*]} (No.~157) or \textbf{[2,3,3,5,5*]}
(No.~180).

\vspace{0.1cm}

\textbf{Degree 3:} Let $Y=\mathbf{X_{141}}, \mathbf{X_{222}}$ or
$\mathbf{X_{411}}$ as defined in Section \ref{s:degree 4}. After
deflation, any base change of degree 3 gives a K3 surface. This is
extremal if and only if the base change is only ramified at the
three cusps (i.e.~the cusps have 5 pre-images in total). We will
refer to the possible base changes as triple covers without
specifying the particular one. They are all defined over $\Q$.

\vspace{0.1cm}

\textbf{Degree 4:} Consider $\mathbf{X_{431}}$ as introduced in
Section \ref{s:degree 6}. We want to apply a base change of degree
4. This gives an extremal K3 fibration, if we select the
ramification index (3,1) at the cusp of the $IV^{\ast}$-fibre and
minimize the number of pre-images of the other two cusps at 4. In
fact, we can adequately choose both base changes $\pi_2$ and
$\pi_4$ from the third section after exchanging cusps. The third
useful base change has ramification index (3,1) at every cusp. It
can be given by
\[
\mathbf{\pi_3}((s:t))=(s^3(s-2t):t^3(t-2s)).
\]
This base change, for instance, gives rise to the constellation
\textbf{[1,3,3,9,IV*]} (No.~233).

\vspace{0.1cm}

\textbf{Degree 5:} For the base changes of the next two
paragraphs, the basic rational elliptic surface will be
$\mathbf{X_{321}}$ (defined in Section \ref{s:inflating}). A base
change of degree 5 with ramification index (2,2,1) at 1 (the cusp
of $III^{\ast}$) leads to a K3 surface. After deflation, only one
$III^{\ast}$ remains in the pull-back. The resulting K3 fibrations
is extremal if and only if the other two cusps have the minimal
number of 4 pre-images. There are five such base changes, all but
one defined over $\Q$. We will only go into detail for one of them
and then list the others:

\begin{itemize}

\item The first base change can be given by
\[\mathbf{\pi_E}((s:t))=(s(s^2-5st+5t^2)^2:4t^5),\] since
$s(s^2-5st+5t^2)^2-4t^5=(s-4t)(s^2-3ts+t^2)^2$. As pull-back we
realize \textbf{[2,4,4,5,III*]} (No.~259) and
\textbf{[1,2,2,10,III*]} (No.~275). \vspace{0.08cm} \item
$\mathbf{\pi_F}((s:t))=(4s^3(3s-5t)^2:t^4(15s+2t))$.
\vspace{0.08cm} \item
$\mathbf{\pi_G}((s:t))=(s^3(4s-5t)^2:t^3(4t-5s)^2)$.
\vspace{0.08cm} \item
$\mathbf{\pi_H}((s:t))=(64t^5:(t-s)^3(9s^2-33st+64t^2))$.
\vspace{0.08cm} \item
$\mathbf{\pi_I}((s:t))=(s^4(s-5t):t^4(2i-11)(5s+(3+4i)t)$ with
$i^2=-1$.

\end{itemize}

\vspace{0.1cm}

\textbf{Degree 6:} Consider a base change of degree 6 with
ramification index (2,2,2) at 1. We apply this base change to
$\mathbf{X_{321}}$. After deflation and transfer of *, we obtain
an elliptic K3 surface whose singular fibres sit above the
pre-images of the other two cusps 0 and $\infty$. Restricting
their number to the minimum 5, we achieve an extremal K3
fibration. The transfer of
* turns a distinct semi-stable fibre $I_n$ into its
non-reduced relative of type $I_n^{\ast}$. Choosing different
$I_n$, one base change gives rise to at most 9 different
fibrations. There are six base changes with the above properties.
We give the three maps without factorization:

\begin{itemize}

\item $\mathbf{\pi_A}((s:t))=(4(s^2-4st+t^2)^3:27t^4s(s-4t))$ or alternatively\newline
$\mathbf{\pi_A'}((s:t))=(4s^3(s-2t)^3:t^4(3s^2-6st-t^2))$.
\vspace{0.08cm}
\item $\mathbf{\pi_B}((s:t))=(-4t^5(6s+t):s^3(2s-5t)^2(s-4t))$.
\vspace{0.08cm}
\item $\mathbf{\pi_C}((s:t))=(-s^4(s^2+2st+5t^2):4t^5(t-2s))$.

\end{itemize}

Let us mention one particular example in more detail: We want to
realize the configuration \textbf{[1,2,3,10,2*]} (No.~148). This
can be achieved as pull-back from $\mathbf{X_{321}}$ via the base
change $\pi_B$. Here we need ramification index (3,2,1) at the
$I_1$-fibre and (5,1) at the $I_2$. The remarkable point about
this construction is that we still have a choice of where to move
the * after the pull-back: We can transfer it either to the $I_2$
at $5/2$ which sits above the original $I_1$, or to the $I_2$ over
$-1/6$ which comes from the $I_2$-fibre of $\mathbf{X_{321}}$.
Using lattice theory, one can easily show that the resulting two
complex surfaces are not isomorphic. This is the only example
where there is such an ambiguity concerning the transfer of *.

We are now in the position to compute all the remaining extremal
elliptic K3 fibrations which can be derived from rational elliptic
surfaces by our simple methods. The following table collects their
configurations together with the number at which they appear in
\cite{SZ}. We further add the Mordell-Weil group $MW$ and the
reduced coefficients of the intersection form $\begin{pmatrix} a &
b\\b & c\end{pmatrix}$ on the transcendental lattice which
determine the isomorphism class of the surface (up to orientation
if there is ambiguity in the sign of $b$). The right-hand part of
the table gives a very brief description of the construction and
the field of definition for the fibration. For shortness we will
indicate the occurrence of a transfer of * only by a * in the
description of the construction. We will not mention deflation.

\vspace{0.2cm}

\begin{table}[ht!]
\begin{footnotesize}
\begin{tabular}{|c||c|c|ccc||c|c|}
\hline No. & Config. & $MW$ & $a$ &  $b$ &  $c$ & Construction & def. \\
\hline
\hline 114 & 1,4,6,6,1* & $\Z/2$ & 12 & 0 & 12 & pull-back from $X_{321}$ via $\pi_A$ * & $\Q$ \\
\hline 115 & 1,5,5,6,1* & $0$ & 20 & 0 & 30 & double cover of [1,3,5,III] * & $\Q$ \\
\hline 116 & 2,4,5,6,1* & $\Z/2$ & 12 & 0 & 20 & pull-back from $X_{321}$ via $\pi_B$ * & $\Q$  \\
\hline 117 & 1,2,7,7,1* & $0$ & 14 & 0 & 28 & double cover of [1,1,7,III] * & $\Q(\sqrt{-7})$ \\
\hline 122 & 2,3,4,8,1* & $\Z/4$ & 6 & 0 & 8 & triple cover of $X_{141}$ & $\Q$ \\
\hline 123 & 2,2,5,8,1* & $\Z/2$ & 8 & 0 & 20 & pull-back from $X_{321}$ via $\pi_C$ * & $\Q$ \\
\hline 127 & 1,3,3,10,1* & $0$ & 6 & 0 & 60 & double cover of [1,3,5,III] * & $\Q$ \\
\hline 128 & 2,2,3,10,1* & $\Z/2$ & 2 & 0 & 60 & pull-back of $X_{321}$ via $\pi_B$ * & $\Q$ \\
\hline 129 & 1,2,4,10,1* & $\Z/2$ & 8 & 4 & 12 & pull-back from $X_{321}$ via $\pi_C$ * & $\Q(\sqrt{-1})$ \\
\hline 132 & 1,2,2,12,1* & $\Z/4$ & 2 & 0 & 6 & triple cover of $X_{141}$ & $\Q$ \\
\hline 133 & 1,1,3,12,1* & $\Z/2$ & 6 & 0 & 6 & triple cover of $X_{141}$ & $\Q$ \\
\hline 135 & 1,1,1,14,1* & $0$ & 6 & 2 & 10 & double cover of [1,1,7,III] * & $\Q(\sqrt{-7})$ \\
\hline 138 & 3,3,5,5,2* & $0$ & 30 & 0 & 30 & double cover of [1,3,5,III] * & $\Q$ \\
\hline 139 & 2,2,6,6,2* & $\Z/2\times\Z/2$ & 6 & 0 & 6 & triple cover of $X_{222}$ & $\Q$ \\
\hline 140 & 2,4,4,6,2* & $\Z/2\times\Z/2$ & 4 & 0 & 12 & triple cover of $X_{222}$ & $\Q$ \\
\hline
\end{tabular}
\end{footnotesize}
\caption{The extremal K3 fibrations with five cusps}
\label{T:extr-5}
\end{table}

\begin{table}%[ht!]
\begin{footnotesize}
\begin{tabular}{|c||c|c|ccc||c|c|}
\hline No. & Config. & $MW$ & $a$ & $b$ & $c$ & Construction & def. \\
\hline
\hline 141 & 1,4,5,6,2* & $\Z/2$ & 4 & 0 & 30 & pull-back from $X_{321}$ via $\pi_B$ * & $\Q$ \\

\hline 142 & 1,1,7,7,2* & $0$ & 14 & 0 & 14 & double cover of [1,1,7,III] * & $\Q(\sqrt{-7})$ \\
\hline 144 & 2,3,3,8,2* & $\Z/2$ & 6 & 0 & 24 & pull-back from $X_{321}$ via $\pi_A$ * & $\Q$ \\
\hline 145 & 1,3,4,8,2* & $\Z/2$ & 4 & 0 & 24 & triple cover of $X_{141}$ * & $\Q$ \\
\hline 146 & 1,2,5,8,2* & $\Z/2$ & 6 & 2 & 14 & pull-back from $X_{321}$ via $\pi_C$ * & $\Q(\sqrt{-1})$ \\
\hline 148 & 1,2,3,10,2* & $\Z/2$ & $\begin{matrix} 6\\4\end{matrix}$ & $\begin{matrix} 0\\2\end{matrix}$ & $\begin{matrix} 10\\16\end{matrix}$ & pull-back from $X_{321}$ via $\pi_B$ * & $\Q$ \\
\hline 149 & 1,1,4,10,2* & $\Z/2$ & 4 & 0 & 10 & pull-back from $X_{321}$ via $\pi_C$ * & $\Q$ \\
\hline 151 & 1,1,2,12,2* & $\Z/2$ & 4 & 0 & 6 & triple cover of $X_{141}$ * & $\Q$ \\
\hline 156 & 3,4,4,4,3* & $\Z/4$ & 8 & 4 & 8 & triple cover of $X_{141}$ & $\Q$ \\
\hline 157 & 2,3,5,5,3* & $0$ & 10 & 0 & 60 & double cover of [1,3,5,III] * & $\Q$ \\
\hline 161 & 2,2,3,8,3* & $\Z/2$ & 4 & 0 & 24 & pull-back from $X_{321}$ via $\pi_A'$ * & $\Q$ \\
\hline 162 & 1,2,4,8,3* & $\Z/4$ & 2 & 0 & 8 & triple cover of $X_{141}$ & $\Q$ \\
\hline 163 & 1,2,2,10,3* & $\Z/2$ & 4 & 0 & 10 & pull-back from $X_{321}$ via $\pi_B$ * & $\Q$ \\
\hline 164 & 1,1,3,10,3* & $0$ & 2 & 0 & 60 & double cover of [1,3,5,III] * & $\Q$ \\
\hline 166 & 1,1,1,12,3* & $\Z/4$ & 2 & 1 & 2 & triple cover of $X_{141}$ & $\Q$ \\
\hline 170 & 3,3,4,4,4* & $\Z/2$ & 12 & 0 & 12 & double cover of [2,3,4,III] * & $\Q$ \\
\hline 171 & 1,1,6,6,4* & $\Z/2$ & 6 & 0 & 6 & pull-back from $X_{321}$ via $\pi_A$ * & $\Q$ \\
\hline 172 & 2,2,4,6,4* & $\Z/2\times\Z/2$ & 2 & 0 & 12 & triple cover of $X_{222}$ * & $\Q$ \\
\hline 173 & 1,2,5,6,4* & $\Z/2$ & 2 & 0 & 30 & pull-back from $X_{321}$ via $\pi_B$ * & $\Q$ \\
\hline 175 & 1,2,3,8,4* & $\Z/2$ & 2 & 0 & 24 & triple cover of $X_{411}$ & $\Q$ \\
\hline 176 & 1,1,2,10,4* & $\Z/2$ & 2 & 0 & 10 & pull-back from $X_{321}$ via $\pi_C$ * & $\Q$ \\
\hline 180 & 2,3,3,5,5* & $0$ & 12 & 0 & 30 & double cover of [1,3,5,III] * & $\Q$ \\
\hline 181 & 1,2,4,6,5* & $\Z/2$ & 4 & 0 & 12 & pull-back from $X_{321}$ via $\pi_B$ * & $\Q$ \\
\hline 182 & 1,1,5,6,5* & $0$ & 4 & 0 & 30 & double cover of [1,3,5,III] * & $\Q$ \\
\hline 184 & 1,2,2,8,5* & $\Z/2$ & 4 & 0 & 8 & pull-back from $X_{321}$ via $\pi_C$ * & $\Q$ \\
\hline 188 & 2,2,4,4,6* & $\Z/2\times\Z/2$ & 4 & 0 & 4 & triple cover of $X_{222}$ & $\Q$ \\
\hline 189 & 1,1,5,5,6* & $0$ & 10 & 0 & 10 & double cover of [1,3,5,III] * & $\Q$ \\
\hline 190 & 1,2,4,5,6* & $\Z/2$ & 2 & 0 & 20 & pull-back from $X_{321}$ via $\pi_B$ * & $\Q$ \\
\hline 191 & 2,2,2,6,6* & $\Z/2\times\Z/2$ & 4 & 2 & 4 & triple cover of $X_{222}$ & $\Q$ \\
\hline 192 & 1,1,4,6,6* & $\Z/2$ & 2 & 0 & 12 & pull-back from $X_{321}$ via $\pi_A'$ * & $\Q$ \\
\hline 201 & 1,1,2,7,7* & $0$ & 6 & 2 & 10 & double cover of [1,1,7,III] * & $\Q(\sqrt{-7})$ \\
\hline 202 & 2,2,3,3,8* & $\Z/2$ & 6 & 0 & 6 & pull-back from $X_{321}$ via $\pi_A$ * & $\Q$ \\
\hline 203 & 1,2,3,4,8* & $\Z/2$ & 4 & 0 & 6 & triple cover of $X_{411}$ * & $\Q$ \\
\hline 204 & 1,2,2,5,8* & $\Z/2$ & 2 & 0 & 4 & pull-back from $X_{321}$ via $\pi_C$ * & $\Q$ \\
\hline 210 & 1,1,3,3,10* & $0$ & 6 & 0 & 6 & double cover of [1,3,5,III] * & $\Q$ \\
\hline 211 & 1,2,2,3,10* & $\Z/2$ & 2 & 0 & 6 & pull-back from $X_{321}$ via $\pi_B$ * & $\Q$ \\
\hline 212 & 1,1,2,4,10* & $\Z/2$ & 2 & 0 & 4 & pull-back from $X_{321}$ via $\pi_C$ * & $\Q$ \\
\hline 215 & 1,1,2,2,12* & $\Z/2$ & 2 & 0 & 2 & triple cover of $X_{411}$ & $\Q$ \\
\hline 216 & 1,1,1,3,12* & $\Z/2$ & 2 & 1 & 2 & triple cover of $X_{411}$ & $\Q$ \\
\hline 218 & 1,1,1,1,14* & $0$ & 2 & 0 & 2 & double cover of [1,1,7,III] * & $\Q$ \\
\hline 223 & 1,3,6,6,IV* & $\Z/3$ & 6 & 0 & 6 & pull-back from $X_{431}$ via $\pi_2$ & $\Q$ \\
\hline 224 & 3,3,4,6,IV* & $\Z/3$ & 6 & 0 & 12 & pull-back from $X_{431}$ via $\pi_4$ & $\Q$ \\
\hline 233 & 1,3,3,9,IV* & $\Z/3$ & 6 & 3 & 6 & pull-back from $X_{431}$ via $\pi_3$ & $\Q$ \\
\hline 234 & 2,2,3,9,IV* & $\Z/3$ & 2 & 0 & 18 & pull-back from $X_{431}$ via $\pi_2$ & $\Q$ \\
\hline 241 & 1,1,2,12,IV* & $\Z/3$ & 2 & 0 & 4 & pull-back from $X_{431}$ via $\pi_4$ & $\Q$ \\
\hline 259 & 2,4,4,5,III* & $\Z/2$ & 4 & 0 & 20 & pull-back from $X_{321}$ via $\pi_E$ & $\Q$ \\
\hline 261 & 1,4,4,6,III* & $\Z/2$ & 4 & 0 & 12 & pull-back from $X_{321}$ via $\pi_F$ & $\Q$ \\
\hline 262 & 2,3,4,6,III* & $\Z/2$ & 6 & 0 & 12 & pull-back from $X_{321}$ via $\pi_G$ & $\Q$ \\
\hline 263 & 2,2,5,6,III* & $\Z/2$ & 8 & 2 & 8 & pull-back from $X_{321}$ via $\pi_H$ & $\Q$ \\
\hline 270 & 2,2,3,8,III* & $\Z/2$ & 2 & 0 & 24 & pull-back from $X_{321}$ via $\pi_F$ & $\Q$ \\
\hline 271 & 1,2,4,8,III* & $\Z/2$ & 4 & 0 & 8 & pull-back from $X_{321}$ via $\pi_I$ & $\Q(\sqrt{-1})$ \\
\hline 275 & 1,2,2,10,III* & $\Z/2$ & 2 & 0 & 10 & pull-back from $X_{321}$ via $\pi_E$ & $\Q$ \\
\hline 276 & 1,1,3,10,III* & $\Z/2$ & 4 & 1 & 4 & pull-back from $X_{321}$ via $\pi_H$ & $\Q$ \\
\hline 298 & 3,3,4,4,II* & $0$ & 12 & 0 & 12 & double cover of [3,4,II,III] * & $\Q$ \\
\hline 299 & 2,2,5,5,II* & $0$ & 10 & 0 & 10 & double cover of [2,5,II,III] * & $\Q$ \\
\hline 301 & 1,1,6,6,II* & $0$ & 6 & 0 & 6 & double cover of [1,6,II,III] * & $\Q$ \\
\hline
\end{tabular}
\end{footnotesize}
\vspace{0.5cm}
\begin{center}
{\sf Tab.~2:} The extremal K3 fibrations with five cusps
\end{center}
\end{table}

Note that No.\,135 and 201 are fibrations on the same complex
surface as No.\,282 and the semi-stable fibration with
configuration [1,1,2,2,4,14] from Section \ref{s:inflating}. Since
this surface cannot be defined over $\Q$ by \cite[Prop.~13.1]{S},
the field of definition $\Q(\sqrt{-7})$ is minimal.

Table \ref{T:extr-5} completes the treatment of extremal elliptic
K3 fibrations which can be derived from rational elliptic surfaces
by direct manipulation of the Weierstrass equation or as pull-back
via a non-general base change. We would like to finish this
section with the following remark. It concerns K3 surfaces which
possess an extremal elliptic fibration with non-trivial
Mordell-Weil group. For every such surface, this paper (combined
with \cite{TY}) gives at least one explicit extremal fibration
which is obtained as pull-back from a rational elliptic surface.
This result might be compared to the idea of elementary fibrations
proposed in \cite[section 6]{MP3}. We should, however, point out
that our pull-backs can in general not be called elementary in the
strict sense of \cite{P},\cite{MP3}.

\vspace{0.8cm}

\textbf{Acknowledgement:} I am indepted to K.~Hulek for his
continuous interest and support. The paper benefitted greatly from
discussions with B.~van Geemen and R.~Kloosterman. I would further
like to thank N.~Yui for drawing my attention to this problem and
H.-C.~Graf von Bothmer for kindly running the Macaulay program.

This paper was partially supported by the DFG-Schwerpunkt 1094
"Globale Methoden in der komplexen Geometrie". The final revision
took place while I enjoyed the hospitality of the Dipartimento di
Matematica "Frederico Enriques" of Milano University. Funding from
the network Arithmetic Algebraic Geometry, a Marie Curie Research
Training Network, is gratefully acknowledged.

\vspace{1cm}

Matthias Sch\"utt\\ Institut f\"ur Algebraische Geometrie\\
Universit\"at Hannover\newline Welfengarten 1\\ 30167 Hannover\\
Germany\\
{\tt schuett@math.uni-hannover.de}

\end{document}